\documentclass{article}
\usepackage{authblk}

\usepackage{times}
\usepackage{amsmath}    
\usepackage{amsfonts,bbm,bm}   
\usepackage{amssymb}
\usepackage{graphicx}
\usepackage{overpic}
\usepackage{caption,subcaption}
\usepackage{longtable}
\usepackage{setspace}
\usepackage{lineno}
\usepackage{color}
\usepackage[utf8]{inputenc}
\usepackage{xspace}

\usepackage{microtype}     
\usepackage{hyperref}      

\usepackage{amsfonts}
\usepackage{amssymb}
\usepackage{multirow}
\usepackage[table,xcdraw]{xcolor}
\usepackage{eurosym}
\usepackage[final]{changes}
\newcommand{\bfh}{{\mathbf h}} 

\newcommand{\bfs}{{\bf s}}

\newcommand{\bfz}{{\mathbf{\tilde{s}}}}
\newcommand{\bfr}{{\mathbf r}}
\newcommand{\rr}{\| \bfr \|}

\newcommand{\st}{\bfs,t}

\newcommand{\sti}[1]{\bfs_{#1}, t_{#1}}

\newcommand{\EE}[1]{\mathcal{E}\left[{\, #1} \, \right]}

\newcommand{\om}{\omega}

\newcommand{\rmx}{{\rm x}}

\newcommand{\HH}{{\mathcal H}}
\newcommand{\E}{\mathrm{e}}
\newcommand{\de}{\delta}

\newcommand{\Rd}{\mathbb{R}^{d}}

\newcommand{\Ha}{\mathcal{H}}
\newcommand{\R}{\mathbb{R}}
\newcommand{\Or}{\mathcal{O}}

\newcommand{\mx}{{m_{\rmx}}}
\newcommand{\bfxv}{\mathbf{x}}

\newcommand{\noi}{\noindent}
\newcommand{\la}{\lambda}
\newcommand{\Samp}{\mathbb{S}}

\newcommand{\beq}{\begin{equation}}
\newcommand{\eeq}{\end{equation}}

\newcommand{\stt}{\textcolor{blue}{ST}\xspace}
\newcommand{\sli}{\textcolor{blue}{SLI}\xspace}

\newcommand{\bmthe}{{\bm \theta}}

\begin{document}

\title{Stochastic Local Interaction Model with Sparse Precision Matrix for  Space-Time Interpolation}

\author{Dionissios~T.~Hristopulos\thanks{\texttt{dionisi@mred.tuc.gr};  Corresponding author}\hspace{2pt}}
\author{Vasiliki D. Agou\thanks{\texttt{vagou@isc.tuc.gr}}}

\affil{Geostatistics Laboratory, School of Mineral Resources Engineering, \\Technical University of Crete, Chania, 73100 Greece}

\maketitle
\begin{abstract}
The application of geostatistical and machine learning methods based on Gaussian processes to big space-time data is beset by the requirement for storing and numerically inverting large and dense covariance matrices. Computationally efficient representations of space-time correlations can be constructed using local models of conditional dependence which can reduce the computational load. We formulate a stochastic local interaction  model for regular and scattered space-time data that incorporates interactions within controlled space-time neighborhoods. The strength of the interaction and the size of the neighborhood are defined by means of kernel functions and adaptive local bandwidths.  Compactly supported kernels lead to finite-size local neighborhoods and consequently to sparse precision matrices that admit explicit expression. Hence, the stochastic local interaction model's requirements for storage are modest and the costly covariance matrix inversion is not needed. We also derive a semi-explicit prediction equation and express the conditional variance of the prediction in terms of the diagonal of the precision matrix. For data on regular space-time lattices,  the stochastic local interaction model is equivalent to a Gaussian Markov Random Field.
\end{abstract}

\newpage

\section{Introduction}

Space-time (\stt) data are becoming available in overwhelming volumes and diverse forms due to the continuing growth of
remote-sensing capabilities, the deployment of low-cost, ground-based sensor networks, as well as the increasing
usage of sensors based on unmanned aerial vehicles, and crowdsourcing~\cite{National13}. The ongoing data explosion has an impact in various fields of science and engineering. The modeling and processing of massive \stt datasets poses
conceptual, methodological, and technical challenges. Sufficiently flexible and computationally powerful
solutions  are not widely available to date, because most existing methods
are not designed for global, high-volume, hyper-dimensional, heterogeneous and uncertain \stt
data. For example, classical geostatistical and machine learning methods~\cite{Chiles12,Rasmussen06} are limited by the cubic dependence of the computational time on data size, which is prohibitive even for large purely spatial data.

The modeling and processing of \stt data require  more advanced methods and computational resources than those that are adequate for  purely spatial data. For example, theories  that simply extend spatial statistics by adding a separable time dimension
are often inadequate for capturing realistic correlations and for analyzing massive \stt data~\cite{Christakos92,Cressie11}. Various methods have been proposed for developing non-separable covariance models~\cite{Deiaco02,Kolovos04,dth17b}. Current methods, whether they are based on geostatistics~\cite{Chiles12}, spatio-temporal statistics~\cite{Gneiting06,Cressie11}, or machine learning~\cite{Rasmussen06} face serious scalability problems. A prevailing obstacle in the processing chain is the computationally demanding iterated inversion of \emph{large covariance (Gram) matrices}~\cite{Rasmussen06,Sun12}. Hence, classical methods executed on standard desktop computers are limited to  datasets with size $N \sim \Or(10^3)-\Or(10^{4})$.  Various approaches for alleviating the dimensionality problem (covariance tapering, composite likelihood, low-rank computations,  stochastic partial differential equation  representation, etc.) have been proposed and developed~\cite{Sun12}.

 In the case of continuum random fields, Gaussian field theories of statistical physics provide models with local structure which is derived from the derivatives of the field~\cite{Mussardo10}. Gaussian Markov random fields (GMRFs) also share the local property, since they are defined in terms of interactions that involve local neighborhoods~\cite{Rue05}. Stochastic local interaction (SLI) models are inspired from GMRFs and Gaussian field theory. They are based on the idea that correlations are generated by interactions between neighboring sites and times. The interactions are incorporated in a precision matrix with simple parametric dependence.

 We present a theoretical framework for the analysis of \stt data that is based on \emph{stochastic local interaction} (\sli) models~\cite{dth15,dth17}. This formulation is useful for filling gaps by interpolation in \stt datasets of environmental  importance. For example,  gaps in records of meteorological
 variables need to be reconstructed  for the evaluation of renewable energy potential at candidate sites~\cite{koutroulis10}, while ground-based rainfall gauge networks often have missing data~\cite{Bardossy14}.  The main idea in \sli is that the \stt correlations are determined by means of \emph{sparse precision matrices} that only involve couplings between near neighbors (in the \stt domain). In contrast with GMRF models that are typically defined on regular lattice data and field theories which are defined on continuum spaces, the \sli framework is suitable for direct application to scattered data and stochastic graph processes.  However, it is also applicable to data on regular space-time lattices.

The remainder of this manuscript is structured as follows: In Section~\ref{sec:sli-model} we define the \stt-\sli model and discuss  its properties. In Section~\ref{sec:prediction} we formulate \stt prediction based on the \sli model, and in Section~\ref{sec:estimation} we discuss parameter estimation from \stt data. Following the theoretical formulation, Section~\ref{sec:data} presents an application of the \sli method to three  datasets which involve simulated \stt data, reanalysis temperature data, and atmospheric ozone measurements. Finally, we present our conclusions and a brief discussion in Section~\ref{sec:conclusions}.

\section{\stt Model based on Stochastic Local Interactions}
\label{sec:sli-model}

A \emph{space-time scalar random field (STRF)} $X(\st;\om) \in \R$
where  $\st \in  \Rd \times \R$ and  $\om \in \Omega $ is defined as a
mapping from the probability space $(\Omega,A,P)$ into the space of
real numbers $\R$. For each \stt coordinate $(\st)$, $ X(\st;\om)$
is a measurable function of $\om$, where $\om$ is the state index~\cite{Christakos92}.
The states (realizations) of the random field $X(\st;\om)$ are real-valued functions $x(\st)$ obtained for a specific $\om$. In the following, the state index $\omega$ is dropped to simplify notation.

We focus on partially sampled realizations $\bfxv = \left(x_{1}, \ldots, x_{N}\right)^\top$ of the random field, where $N \in \mathbb{N}$ is the sample size. The vector $\bfxv$ comprises the
field values at the \stt point set $\Samp= \{ (\sti{1}), \ldots, (\sti{N}) \}$.
The point set is assumed to be quite general; it may represent a time sequence of lattice sites, randomly scattered points in space and time, or a collection of time series at random locations in space.

\subsection{Energy of the exponential joint density}
The \sli  model is based on a joint pdf defined by the \emph{Boltzmann-Gibbs} exponential distribution
\beq
\label{eq:gibbs-pdf}
f_{\rmx}(\bfxv;\bmthe) = \frac{\E^{-\Ha(\bfxv; \bmthe)}}{Z(\bmthe)},
\eeq

where $\Ha(\cdot; \cdot)$ is an energy function that represents the ``cost'' of a specific configuration, $\bmthe$ is a vector of model parameters, and $Z(\bmthe)$ is the normalizing factor known as partition function.

The energy-based approach  is commonly used in statistical physics~\cite{Kardar07,Mussardo10}. Its main advantage is that it expresses statistical dependence in terms of interactions between space locations and time instants which can be local, without recourse to the concept of the covariance matrix. Depending on the form of the interactions involved in the energy, both Gaussian and non-Gaussian probability density functions can be obtained. The most famous example of non-Gaussian dependence is the magnetic Ising model~\cite{Ising25} which was introduced in spatial statistics by Besag~\cite{Besag74}. While non-Gaussian models are definitely interesting, their Gaussian counterparts lead to explicit predictive expressions and uncertainty estimates based on the conditional variance. Hence, herein we focus on a Gaussian \sli model.

We assume that  $\Ha(\bfxv;\bmthe)$ satisfies the following properties  for any vector $\bfxv \in \mathbb{R}^{N}$ and $N \in \mathbb{N}$:

\begin{enumerate}
\item \emph{Gaussianity:} $\Ha(\bfxv;\bmthe)$ is a quadratic function of the data vector $\bfxv$ that can be expressed as
\beq
\label{eq:Gaussian-energy}
\Ha(\bfxv;\bmthe) = \frac{1}{2} (\bfxv - \mathbf{\mx})^{\top} \mathbf{J}(\bmthe') (\bfxv - \mathbf{\mx}),
\eeq
where $\mathbf{\mx}= \left( \mx_{;1}, \ldots, \mx_{;N} \right)^\top$ is a vector of mean (trend) values such that $\mx_{;,i}= \EE{X(\bfs_{i},t_{i})}$, where $\EE{\cdot}$ is the expectation operator.
 On the other hand, $\mathbf{J}(\bmthe')$ is the $N \times N$ precision (or interaction) matrix. The latter depends on the parameter vector $\bmthe' = \bmthe \setminus \{b_{1}, \ldots, b_{K}\}$ which excludes the trend coefficients. The vector $\mathbf{\mx}$  incorporates both periodic and aperiodic trend components.

\item \emph{Positive-definiteness:} $\Ha(\bfxv;\bmthe) >0$ for all $\bfxv$ that are not identically equal to zero.  This is equivalent to the \emph{precision matrix} $\mathbf{J}(\bmthe)$ being a positive-definite matrix.

\item \emph{Sparseness:} $\mathbf{J}(\bmthe')$ is a \emph{sparse matrix} that incorporates the local interactions.
\end{enumerate}

More specifically, we focus on the following  \sli energy function which satisfies the properties of Gaussianity, positive-definiteness and sparseness:
\beq
\label{eq:sli-energy}
\Ha(\bfxv;\bmthe) = \frac{1}{2\la} \left[ \sum_{n=1}^{N} \frac{1}{N}(x_{n} - \mx_{;n})^{2} +  c_{1} \big\langle \left( {x'}_{n} - {x'}_{k}\right)^{2} \big\rangle  \right].
\eeq
We assume that the mean is modeled by means of a trend function which can be expressed as $\mx(\bfs,t) = \sum_{k=1}^{K} b_{k} f_{k}(\bfs, t)$  in terms of a suitable \stt function basis $\{f_{k}(\bfs, t)\}_{k=1}^{K}$, where $\{ b_{k}\}_{k=1}^{K}$ is a set of real-valued trend coefficients and $f_{k}: \Rd\times \R \to \R$, for $k=1, \ldots, K$.

The variables $x_{n}, x_{k}$ stand for $x(\bfs_{n},t_{n})$ and $x(\bfs_{k},t_{k})$ respectively, where $n, k=1, \ldots N$ while ${x'}_{n}, {x'}_{k}$ represent the residuals after the trend values are removed. The term $\big\langle \left( {x'}_{n} - {x'}_{k}\right)^{2} \big\rangle$ represents a weighted average of the squared increments. However, instead of focusing on all $\Or(N^2)$ pairs, the average defined below selects only pairs within a local neighborhood around each point $\bfs_n$.

The \emph{\sli parameter vector} $\bmthe$ includes  the trend coefficients $\{ b_{k}\}_{k=1}^{K}$, the overall scale parameter $\la$ (which is proportional to the variance), and the increment coefficient $c_{1}$  (a dimensionless factor that multiplies the contribution from the squares of the  increments). The vector $\bmthe$ includes additional parameters that determine the local \stt neighborhoods used in the average of the squared increments $\langle \cdot \rangle$. The average is defined in~\eqref{eq:squared-incr} below.

\subsection{Kernel-based averaging}
The weights in the average of the squared increments are defined by means of the \emph{Nadaraya-Watson} equation~\cite{Nadaraya64,Watson64}, i.e.,
\beq
\label{eq:squared-incr}
\langle \left( {x'}_{n} - {x'}_{k}\right)^{2} \rangle =
\frac{\sum_{n=1}^{N}\sum_{k=1}^{N} w_{n,k}\left( {x'}_{n} - {x'}_{k}\right)^{2} }{\sum_{n=1}^{N}\sum_{k=1}^{N} w_{n,k}}.
\eeq

The coefficients $w_{n,k}$ are defined in terms of \emph{compactly supported \stt kernel functions} $K(\cdot,\cdot): \R^{q}\times \R^{q} \to \R$, where $q=d+1$ for an \stt kernel, $q=d$ for a spatial kernel, and $q=1$ for a temporal kernel. Kernel functions are symmetric, real-valued functions; herein they are assumed  to take values in the interval $[0, 1]$ without loss of generality. Moreover, we will assume spatially homogeneous and temporally stationary kernel functions, i.e.,  $K(\bfs_{1},\bfs_{2})=K(\bfs_{1}-\bfs_{2})$, $K(t_{1},t_{2})=K(t_{1}-t_{2})$, and  $K(\bfs_{1},t_{1};\bfs_{2},t_{2})=K(\bfs_{1}-\bfs_{2},t_{1}-t_{2})$. Furthermore, it will be assumed for simplicity that the kernel function depends only on the magnitude of the \stt distance.

\subsection{Definition of space-time distance}
The space-time distance used in the kernel weights determines the structure of correlations that we impose in the space-time domain. Both separable and non-separable space-time metric distances are possible as discussed below.

\emph{Composite space-time distance:} In this case the spatial and temporal coordinates are intertwined in the distance metric. For example, the  differential of the space-time distance between two points using the Riemannian metric is
\beq
dq = \sqrt{\sum_{i=1}^{d+1}\sum_{j=1}^{d+1} g_{i,j} dz^{(i)} dz^{(j)}},
\eeq
where $\{ g_{i,j} \}_{i,j=1}^{d+1}$ are the elements of the metric tensor $\mathbf{g}$, $\{ dz^{(i)} \}_{i=1}^{d}$,  are the differentials of the spatial distance in the $d$ orthogonal directions, and $dz^{(d+1)}$ is the time differential~\cite{awr00,Christakos17}.

\begin{figure}[ht]
  \centering
  \begin{subfigure}[b]{0.46\textwidth}
    \centering\includegraphics[width=\linewidth]{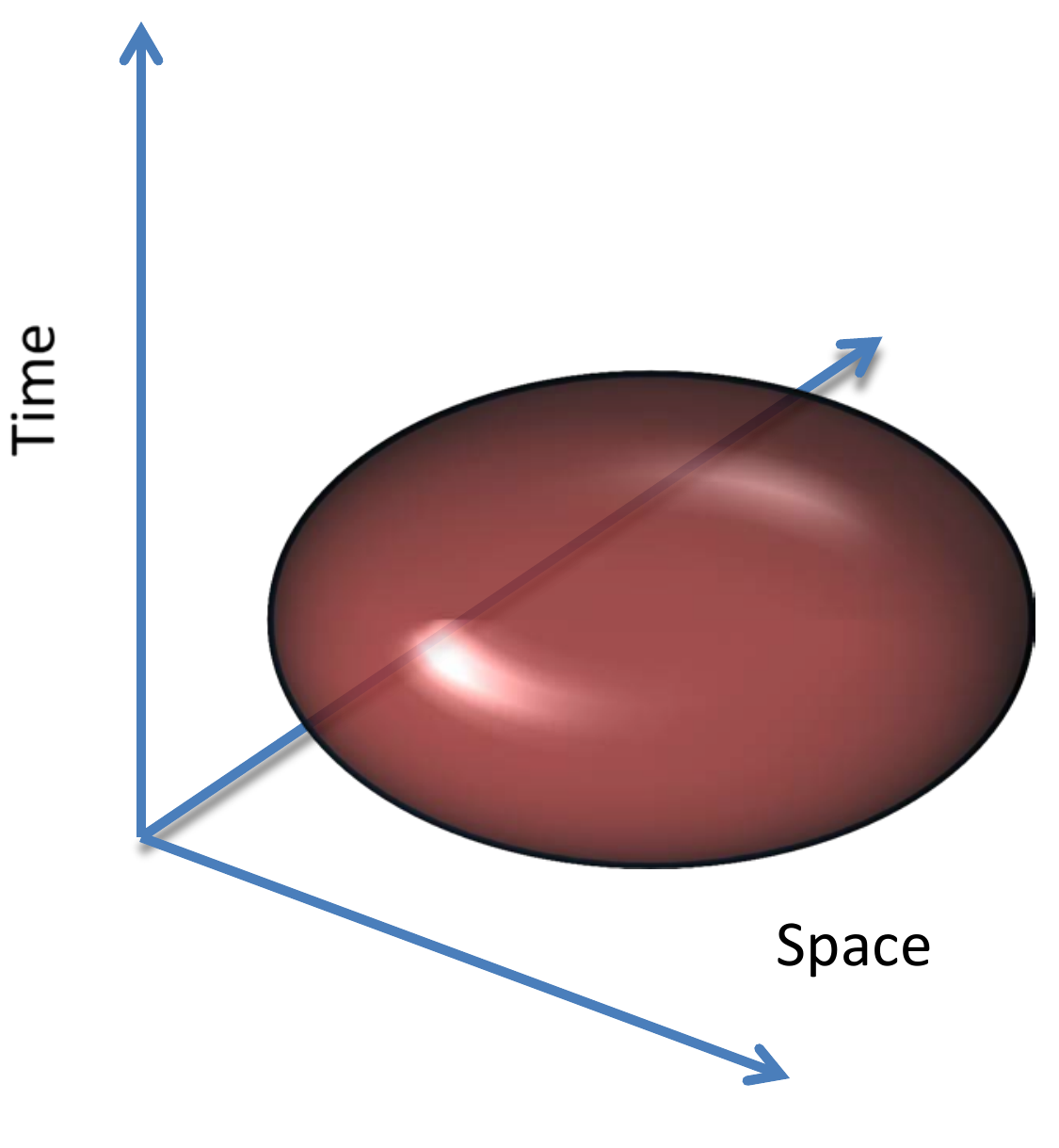}
    \caption{Composite\label{fig:st-kernel-compos}}
  \end{subfigure}%
  \begin{subfigure}[b]{0.46\textwidth}
    \centering\includegraphics[width=\linewidth]{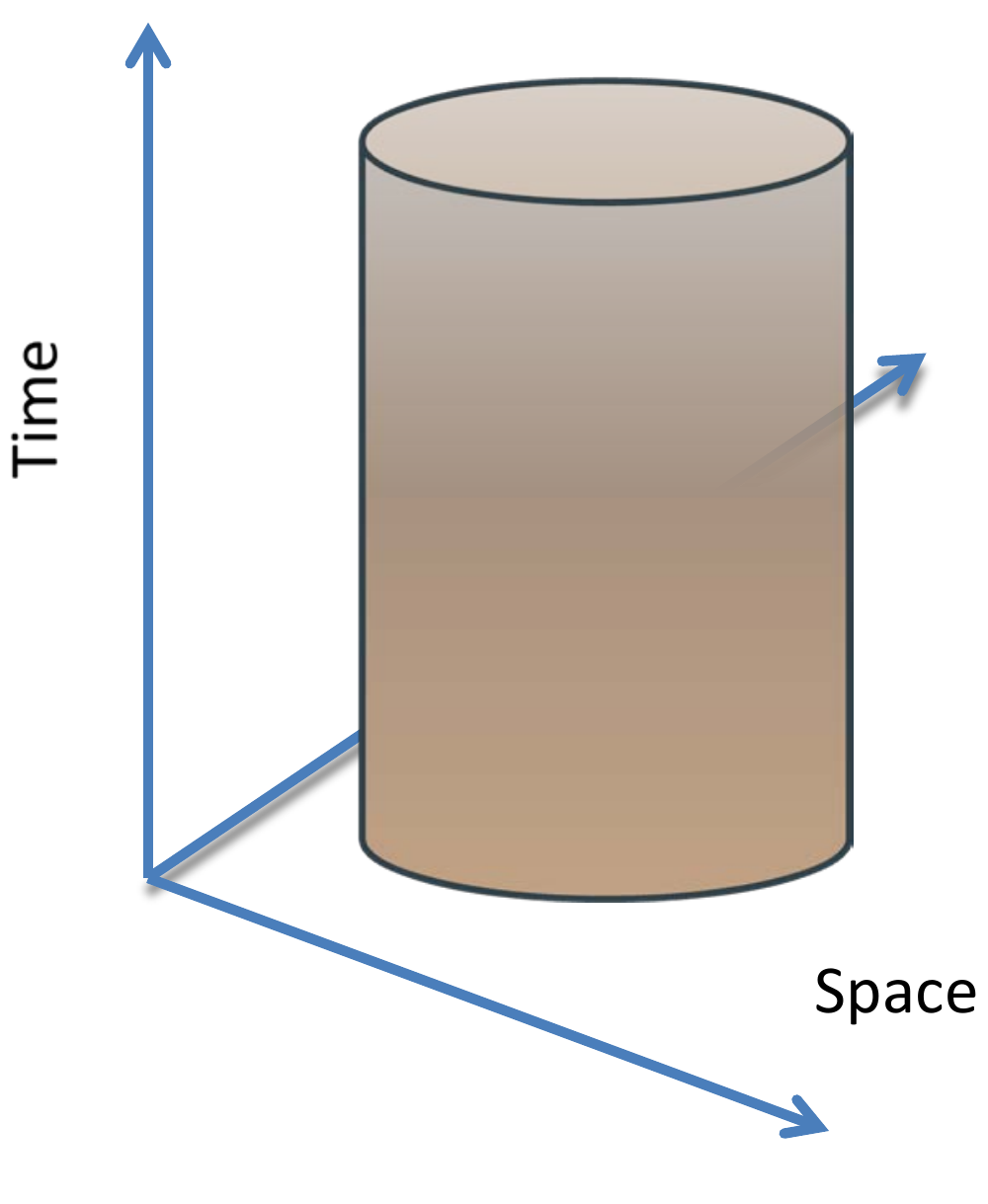}
    \caption{Separable\label{fig:st-kernel-separ}}
  \end{subfigure}
  \caption{Schematics of kernel-based neighborhoods for composite (left) and separable (right) space-time structures.}
\end{figure}

In the Euclidean case the metric tensor $\mathbf{g}$,  is given by
\beq
g_{i,j}=\de_{i,j} \left[ 1 + \left(\alpha -1 \right)\de_{i,d+1}\right], \; i,j=1, \ldots, d+1,
\eeq
where $\alpha$ is a  parameter that controls the contribution of the time lag in the composite distance.
The kernel coefficient based on the composite Euclidean  metric can be expressed as
\beq
\label{eq:k-weight-compos}
w_{n,k} = K\left( \frac{\sqrt{ \bfr_{n,k}^{2} + \alpha^{2}\tau_{n,k}^{2} }}{h_{s,n}}\right).
\eeq
In~\eqref{eq:k-weight-compos}  $\bfr_{n,k} = \left( \bfs_{n} - \bfs_{k}\right)$ is the  spatial lag between the initial point $\bfs_{n}$ and the target point $\bfs_{k}$, and $h_{s,n}$ is the \emph{local spatial bandwidth} at $\bfs_{n}$. In addition, $\tau_{n,k}= t_{n} - t_{k}$ is the temporal lag between the initial and target times. The space-time distance for the composite metric leads to  ellipsoidal  neighborhoods as  shown in the schematic of  Fig.~\ref{fig:st-kernel-compos}. The temporal bandwidth in this case is $h_{t,n}=h_{s,n}/\alpha$.

\emph{Separable space-time distance:} The coefficients $w_{n,k}$ for a  separable space-time neighborhood are defined as
\beq
\label{eq:k-weight-separ}
w_{n,k}= K\left( \frac{\| \bfr_{n,k}\|}{h_{s,n}}\right)\, K\left( \frac{| \tau_{n,k}|}{h_{t,n}}\right), \;n, k=1, \ldots, N.
\eeq
In the weight equation~\eqref{eq:k-weight-separ}   $h_{s,n}$ is the \emph{local spatial bandwidth} at $\bfs_{n}$ and $h_{t,n}$ is the \emph{temporal bandwidth}.
The space-time distance for the separable space-time metric leads to cylindrical  neighborhoods as shown in Fig.~\ref{fig:st-kernel-separ}.

\subsection{Definition of bandwidths}

For each \stt point $\{ (\bfs_{n}, t_{n})\}_{n=1}^{N}$, the \emph{spatial bandwidth} $h_{s,n}$ is determined from the geometry of the sampling network around the spatial point $\bfs_{n}$, while the temporal bandwidth $h_{t,n}$ is based on the time neighborhood around $t_{n}$.  In general, this means that the number of bandwidth parameters scales linearly with the sampling size, leading to an under-determined estimation problem when the additional parameters are accounted for.

To simplify the bandwidth estimation we use a trick that reduces the dimensionality of the problem. We assign to each point a bandwidth which is proportional to the spatial distance $D_{n, [K_{s}]}(\Samp)$ between this point and its $K_{s}$-nearest neighbor in the point set $\Samp$.  Thus, it holds that $h_{s,n}=\mu_{s} \, D_{n, [K_{s}]}(\Samp)$, where typically $K_{s}=2, 3, 4$, and $\mu_{s} >0$ is a dimensionless spatial bandwidth parameter to be estimated from the data.

In the case of the \emph{composite space-time} distance the temporal bandwidths $h_{t,n}$ are determined from the $h_{s,n}$ and the additional parameter $\alpha$.
For a \emph{separable \stt} distance metric, the temporal bandwidths are determined by means of
$h_{t,n}=\mu_{t} \, \tilde{D}_{n, [K_{t}]}(\Samp)$, where $K_{t}$ is the order of the temporal neighbor and $\mu_{t}>0$ is a dimensionless temporal bandwidth parameter.  This definition of the temporal bandwidth in the case of uniform time step implies uniform bandwidths for all except the initial and final times, where the bandwidth is automatically increased to account for the missing left and right neighbors respectively.

\subsection{Properties of kernel weights}

The kernel-average of the squared increments~\eqref{eq:squared-incr} can be expressed in terms of normalized weights $u_{n,k}$ as follows
\begin{subequations}
\label{eq:squared-increments-norm-weights}
\begin{align}
\label{eq:xnxk2-av}
\langle \left( {x'}_{n} - {x'}_{k}\right)^{2} \,\rangle  & =
\sum_{n=1}^{N}\sum_{k=1}^{N} u_{n,k}\left( {x'}_{n} - {x'}_{k}\right)^{2},
\\
\label{eq:unk}
u_{n,k} & = \frac{w_{n,k}}{\sum_{n=1}^{N}\sum_{k=1}^{N} w_{n,k}}.
\end{align}
\end{subequations}

\emph{Normalization:}  The definition~\eqref{eq:unk} of the kernel weights implies that
\[
\sum_{n=1}^{N}\sum_{k=1}^{N} u_{n,k}=1.
\]

\emph{Asymmetry:} The  definition~\eqref{eq:unk} of the bandwidths is based on the local \stt  neighborhood. This implies that the \emph{spatial weights} are in general asymmetric, i.e.,  $w_{n,k} \neq w_{k,n}$ if $\bfs_{n} \neq \bfs_{k}$, since the sampling density around the point $\bfs_{n}$ can be quite different than around the point $\bfs_{k}$.

\emph{Non-separability:} The kernel weights $u_{n,k}$ are non-separable for both the composite and the separable \stt distance metrics. In the first case this  is obvious from the definition~\eqref{eq:k-weight-compos}.
In the second case, even though the $w_{n,k}$ are separable, the normalized weights $u_{n,k}$ are non-separable functions of space and time due to the kernel summation in the denominator of~\eqref{eq:squared-incr}.

\emph{Robustness with respect to general distance metrics:} Regardless of the distance metric used, the kernel-based weights $u_{n,k}$ are non-negative. This implies that the \sli energy function~\eqref{eq:sli-energy} is positive, and consequently the precision matrix is positive definite. Hence, general distance metrics, e.g., Manhattan (also known as city block and taxicab) distance, can be used in the \sli model.

In the following we develop the \sli formalism for a separable space-time metric structure.

\subsection{Squared increments for separable space-time metric}
In this section we formulate the
average squared increments for  separable space-time kernel functions using matrix operations.

First, we define the square kernel matrices $\mathbf{K}_{s}$ of dimension $N_{s} \times N_{s}$ and $\mathbf{K}_{t}$ of dimension $N_{t} \times N_{t}$ as follows

\begin{subequations}
\beq
\mathbf{K}_{s} = \left[     \begin{array}{cccc}
                            K\left(\frac{\|\bfr_{1,1} \|}{h_{s,1}}\right) &   \ldots &  K\left(\frac{\|\bfr_{1,N_{s}} \|}{h_{s,1}}\right) \\
                             \vdots & \vdots  & \vdots \\
                             K\left(\frac{\|\bfr_{N_{s},1} \|}{h_{s,N_{s}}}\right)  & \ldots &  K\left(\frac{\|\bfr_{N_{s},N_{s}} \|}{h_{s,N_{s}}}\right) \\
                         \end{array}
                \right],
\eeq

\beq
 \mathbf{K}_{t} = \left[     \begin{array}{cccc}
                            K\left(\frac{\|\tau_{1,1} \|}{h_{t,1}}\right) &   \ldots &  K\left(\frac{\|\tau_{1,N_{t}} \|}{h_{t,1}}\right) \\
                            \vdots & \vdots   & \vdots \\
                             K\left(\frac{\|\tau_{N_{t},1} \|}{h_{t,N_{t}}}\right) & \ldots &  K\left(\frac{\|\tau_{N_{t},N_{t}} \|}{h_{t,N_{t}}}\right) \\
                         \end{array}
 \right].
\eeq

Then, the $N \times N$ matrix $\mathbf{W}$ of \stt kernel weights is given by the following Kronecker product (denoted by $\otimes$):
\beq
\label{eq:matrix-W}
\mathbf{W} = \mathbf{K}_{s} \otimes \mathbf{K}_{t}.
\eeq
\end{subequations}
For compactly supported kernel functions the matrix $\mathbf{W}$ given by~\eqref{eq:matrix-W} is sparse.

 The matrix $\mathbf{U}$ of the \emph{normalized kernel weights} is then defined by means of
\begin{subequations}
\beq
\label{eq:U}
\mathbf{U} = \frac{\mathbf{W}}{\| \mathbf{W} \|_{1}},
\eeq
where the   denominator $\|\mathbf{W} \|_{1}$ represents the entry-wise $L_{1}$ norm of the matrix $\mathbf{W}$ and  is given by
\beq
\| \mathbf{W} \|_{1}= \sum_{k=1}^{N} \sum_{l=1}^{N} | W_{k,l} |.
\eeq
\end{subequations}

In terms of the above matrices, the average squared increment~\eqref{eq:squared-increments-norm-weights} is expressed as follows

 \beq
 \label{eq:squared-increments-matrix}
\langle \left( {x'}_{n} - {x'}_{k}\right)^{2} \rangle = \left\| \left[ \left( \bfxv' \otimes \mathbf{1} \right) -\left( \bfxv' \otimes \mathbf{1} \right)^\top  \right] \circ \mathbf{U}  \circ \left[\left( \bfxv' \otimes \mathbf{1} \right) -\left( \bfxv' \otimes \mathbf{1} \right)^\top  \right] \right\|_{1}
\eeq

\noi where $\mathbf{1} = (1,  \ldots, 1)^\top$ is the $N \times 1$ vector of ones, and $\circ$ denotes the Hadamard product, i.e., $\left[ \mathbf{A} \circ\mathbf{B} \right]_{i,j} = A_{i,j} B_{i,j}$.

The computational complexity of the operations in~\eqref{eq:squared-increments-matrix}  is $\Or(N^2)$,
if the sparsity of the matrix $\mathbf{W}$ is not taken into account. However, the numerical complexity can be improved using sparse-matrix operations. We have implemented all the calculations which involve the precision matrix using sparse matrix functionality.

\subsection{Precision matrix formulation}
In light of the above definitions, the \sli energy function~\eqref{eq:sli-energy} involves the following parameter vector

\beq
\label{eq:bmthe}
\bmthe = \left( b_{1}, \ldots, b_{K}, \la, c_{1}, \mu_{s}, \mu_{t}, K_{s}, K_{t} \right)^\top,
\eeq
where $\{ b_{k} \}_{k=1}^{K}$ are the coefficients of the trend model, $\lambda$ is the \sli scaling factor, $c_{1}$ is the square increment coefficient, $\mu_{s}, \mu_{t}$ the dimensionless scaling factors used to determine the bandwidths, and $K_{s}, K_{t}$ are the orders of spatial and temporal near neighbors respectively.

The \sli energy function~\eqref{eq:sli-energy} can be transformed into a quadratic energy functional, i.e., of the form of equation~\eqref{eq:Gaussian-energy}, by defining the precision matrix ${\mathbf J}(\bmthe')$ as follows
 \begin{align}
\label{eq:prec-mat} {\mathbf J}(\bmthe')    =     \frac{1}{\lambda } &
\,\left\{ \frac{{\bf I}_{N}}{N} + c_1 \, {\mathbf J}_{1}({\bfh}; \bmthe'') \right\},
\end{align}

\noindent where  ${\bf I}_{N}$ is the $N\times N$ identity matrix: $[{\bf I}_{N}]_{i,j}=1$ if $i=j$ and
$[{\bf I}_{N}]_{i,j}=0$ otherwise. The precision matrix ${\mathbf J}(\bmthe')$ involves the parameter vector $\bmthe'=\left(\la, c_{1}, \mu_{s}, \mu_{t}, K_{s}, K_{t} \right)^\top$. The matrix ${\mathbf J}_{1}({\bfh}; \bmthe'')$ is derived from the average squared increments~\eqref{eq:squared-increments-norm-weights}, and $\bmthe'' = \left(\mu_{s}, \mu_{t}, K_{s}, K_{t} \right)^\top$ is the parameter vector which determines the kernel bandwidths. The precision matrix is thus expressed in terms of the normalized weights $u_{n,k}$  according to

\begin{align}
\label{eq:Jtilde}
[{\mathbf J}_{1}({\bfh}; \bmthe'')]_{n,k} & = - u_{n,k} - u_{k,n} +
 [{\bf I}_{N}]_{n,k} \, \sum_{l=1}^{N} \left( u_{n,l} + u_{l,n}\right),
 \end{align}
 where the normalized weights $u_{n,k}$ are given by~\eqref{eq:unk}. Hence, the precision matrix is determined by the sampling pattern, the kernel functions, and the bandwidths.


\section{\stt Prediction}
\label{sec:prediction}
In this section we consider \stt prediction by means of the \sli model at the set of  space time points $\mathbb{G}=\{ \bfz_{p}\}_{p=1}^{P}$, where $\bfz_{p}=(\tilde{\bfs}_{p}, \tilde{t}_{p})$, assuming that the model parameters are known.  It is further assumed that the sets $\Samp$ and $\mathbb{G}$ are disjoint. For example, the set $\mathbb{G}$ could comprise all the nodes of a regular map grid at a time instant $t_{p}$ for which measurements are not available. Alternatively, $\mathbb{G}$ could comprise all the nodes of an irregular spatial sampling network at a time instant with no measurements.

\subsection{\sli energy function including prediction set}
The \sli energy function that incorporates the prediction sites is given by straightforward extension of~\eqref{eq:sli-energy}. Thus, the following expression that involves block vectors of sampling and prediction sites and respective precision  block matrices is obtained
\beq
\label{eq:sli-energy-pred}
\HH (\bfxv, \bfxv_{\mathbb{G}} ;\bmthe^{\ast})  =\frac{1}{2}   \left[ \begin{array}{cc}
                                {\bfxv'}^\top & {\bfxv}'_{\mathbb{G}}
                                \end{array} \right] \, \left[  \begin{array}{cc}{\mathbf J}_{\Samp,\Samp} & {\mathbf J}_{\Samp,\mathbb{G}} \\
                                {\mathbf J}_{\mathbb{G},\Samp} & {\mathbf J}_{\mathbb{G},\mathbb{G}}
            \end{array}   \right]
\left[ \begin{array}{c}
  {\bfxv'} \\
  {\bfxv}'_{\mathbb{G}}
\end{array}\right],
\eeq
where $\bfxv'=\bfxv - \mathbf{\mx}$ is the detrended data vector, ${\bfxv}'_{\mathbb{G}}= {\bfxv}_{\mathbb{G}}-\mathbf{\mx}$ is the fluctuation vector at the prediction points, and $\bmthe^{\ast}$ is the estimate of the parameter vector based on the data. Let the sets $A, B$ denote either of the disjoint sets $\Samp$ or $\mathbb{G}$. Then, the block precision matrices ${\mathbf J}_{A,B}$  are  expressed as
\beq
\label{eq:sli-J-prediction}
{\mathbf J}_{A,B}({\bmthe'}^\ast) = \frac{1}{\la} \left[ c_{0} \mathbf{I} + c_{1} {\mathbf J}_{A, B}^{(1)}({\bmthe''}^\ast) \right].
\eeq
The block sub-matrices ${\mathbf J}_{A, B}^{(1)}$ are defined as follows:
\begin{subequations}
\label{eq:block-J}
\begin{align}
\left[ {\mathbf J}_{\Samp,\Samp}^{(1)}\right]_{n,k} = & - u_{n,k} - u_{k,n},\, & n, k=1, \ldots, N, \, n \neq k
 \\
 \left[ {\mathbf J}_{\Samp,\Samp}^{(1)}\right]_{n,n} = &  \sum_{l=1\neq n}^{N} \left( u_{n,l} + u_{l,n}\right) + \sum_{p=1}^{P} \left( u_{n,p} + u_{p,n}\right),\, & n=1, \ldots, N,
 \\
 \label{eq:block-J-SG}
 \left[{\mathbf J}_{\Samp,\mathbb{G}}^{(1)}\right]_{n,p} = & - u_{n,p} - u_{p,n}, \,& n=1, \ldots, N, \; p=1, \ldots P, \\[1ex]
 \label{eq:block-J-GS}
 \left[ {\mathbf J}_{\mathbb{G},\Samp}^{(1)}\right]_{\quad}= &  {\mathbf J}_{\Samp,\mathbb{G}}^{(1)\top},
 \\
 \label{eq:block-J-GG-nondiag}
 \left[ {\mathbf J}_{\mathbb{G},\mathbb{G}}^{(1)}\right]_{p,q} = & - u_{p,q} - u_{q,p},\, & p \neq q=1, \ldots, P,
 \\
 \label{eq:block-J-GG-diag}
 \left[ {\mathbf J}_{\mathbb{G},\mathbb{G}}^{(1)}\right]_{p,p}  = & \sum_{l=1}^{N} \left( u_{p,l} + u_{l,p}\right) + \sum_{q\neq p=1}^{P} \left( u_{p,q} + u_{q,p}\right), \, & p=1, \ldots, P.
\end{align}
\end{subequations}

\subsection{Prediction based on stationary point of the energy}
The  Boltzmann-Gibbs pdf of the field at the prediction sites conditional on the data is given by $\exp\left[ -\HH(\bfxv, {\bfxv}_\mathbb{G} ; \bmthe^{\ast})\right]/Z(\bmthe^{\ast})$. The prediction $\hat{\bfxv}_\mathbb{G}$ maximizes the pdf, which is equivalent to minimizing the energy, i.e.,
\beq
\hat{\bfxv}_\mathbb{G} = \arg \min_{{\bfxv}_\mathbb{G}} \HH(\bfxv, {\bfxv}_\mathbb{G} ; \bmthe^{\ast}).
\eeq
The \sli energy~\eqref{eq:sli-energy-pred} can be further expressed in terms of the precision matrix as follows
\[
\HH(\bfxv, {\bfxv}_\mathbb{G} ; \bmthe^{\ast}) = \HH_{s}(\bfxv; \bmthe^{\ast}) +
\frac{1}{2}\left( {\bfxv'}_{\mathbb{G}}^\top {\mathbf J}_{\mathbb{G},\Samp} {\bfxv'} + {\bfxv'}^\top {\mathbf J}_{\Samp,\mathbb{G}}{\bfxv}'_{\mathbb{G}}
+ {\bfxv'}_{\mathbb{G}}^\top {\mathbf J}_{\mathbb{G},\mathbb{G}} {\bfxv'}_{\mathbb{G}}\right),
\]
where  $\HH_{s}(\bfxv; \bmthe^{\ast})= {\bfxv'}^\top{\mathbf J}_{\Samp,\Samp}\, \bfxv'/2$ depends only on the data and is thus irrelevant for the prediction. The condition for a stationary point of  the energy function  is

\beq
\label{eq:sli-energy-stationary}
\frac{\partial \HH(\bfxv, {\bfxv}_\mathbb{G} ; \bmthe^{\ast}) }{\partial {x'}_{p}}=0, \; \mbox{for all} \; \bfz_{p} \in \mathbb{G}.
\eeq

The Hessian  of the energy is $\nabla' \nabla' \HH(\bfxv, {\bfxv}_\mathbb{G} ; \bmthe^{\ast})$, where the prime denotes differentiation with respect to $\bfxv'$.  For the stationary point to represent a minimum of the energy (and thus  a maximum of the Boltzmann-Gibbs pdf), $\nabla' \nabla' \HH(\bfxv, {\bfxv}_\mathbb{G} ; \bmthe^{\ast})$ must be positive definite. From~\eqref{eq:sli-energy-pred} it follows that $\nabla' \nabla' \HH(\bfxv, {\bfxv}_\mathbb{G} ; \bmthe^{\ast})= {\mathbf J}_{\mathbb{G},\mathbb{G}}$. Since the \sli precision matrix is positive definite by construction, so is the Hessian as well.

Finally, the \sli prediction is given by the following equation
\beq
\label{eq:sli-prediction}
\hat{\bfxv}_\mathbb{G}(\bmthe^{\ast} | \bfxv)  = \mathbf{\mx} -\,{\mathbf J}_{\mathbb{G},\mathbb{G}}^{-1}({\bmthe'}^\ast) \, {\mathbf J}_{\mathbb{G},\Samp}({\bmthe'}^\ast) \,\bfxv',
\eeq
where $\mathbf{\mx}$ is the $P \times P$ diagonal trend matrix, i.e.,
$[\mathbf{\mx}]_{p,q}=\de_{p,q} \mx(\bfs_{p}, t_{p})$ and the precision matrices ${\mathbf J}_{\mathbb{G},\mathbb{G}}$ and ${\mathbf J}_{\mathbb{G},\Samp}$ are defined by means of~\eqref{eq:sli-J-prediction} and~\eqref{eq:block-J-SG}-\eqref{eq:block-J-GG-diag}.

Note that due to the matrix product ${\mathbf J}_{\mathbb{G},\mathbb{G}}^{-1} \, {\mathbf J}_{\mathbb{G},\Samp}$ and in light of~\eqref{eq:sli-J-prediction} the \sli prediction is independent of the scale parameter $\la$. This property is analogous to the independence of the kriging prediction from the variance, since the latter is proportional to $\la$.

\subsection{Prediction intervals}
Since the precision matrix of the \sli model is known, it is straightforward to obtain the \emph{conditional variance} at the prediction sites using the result known in Markov random field theory~\cite{Rue05}. Hence,

\beq
\label{eq:sli-condi-variance}
\sigma^{2}_{\mathrm{sli}}(\bfz_{p}) = \frac{1}{J_{p,p}(\bmthe^\ast)}, \; \bfz_{p} \in  \mathbb{G},
\eeq
where $J_{p,p}(\bmthe^\ast)$ is the $p$-th diagonal entry of the precision matrix
${\mathbf J}_{\mathbb{G},\mathbb{G}}$ which is determined from~\eqref{eq:sli-J-prediction} and~\eqref{eq:block-J-GG-diag}.

Based on the above,  \emph{prediction intervals} at the site $\bfz_{p} \in \mathbb{G}$   can be constructed as follows
\[
[\hat{x}_{p}- z_{q}\sigma_{\mathrm{sli}}(\bfz_{p}),  \; \hat{x}_{p} + z_{q}\sigma_{\mathrm{sli}}(\bfz_{p})],
\]
where  $0 \le q\le 1$ is a specified level (e.g., $q=0.95$), and $z_{q}$ is the respective standard $z$-score.

\section{Parameter Estimation}
\label{sec:estimation}

We use maximum likelihood estimation (MLE) to estimate the \sli model parameter vector~\eqref{eq:bmthe}.
The orders of the spatial and temporal neighbors $K_s$ and $K_t$ are set in advance to low integer values larger than one. This does not have a serious impact on the results, since the bandwidth parameters $\mu_{s}, \mu_{t}$ compensate for the choice of the neighbor order.

The maximization of the \sli likelihood ${\mathcal{L}}(\bmthe;\bfxv )$ is equivalent to  minimizing the negative log-likelihood (NLL). In light of equations~\eqref{eq:gibbs-pdf} and~\eqref{eq:Gaussian-energy}, the NLL is given by
\begin{equation}
\label{eq:loglik}
-\ln{\mathcal{L}}(\bmthe;\bfxv ) = \Ha(\bfxv ;\bmthe) + \ln{Z(\bmthe)} =
\frac{1}{2} (\bfxv - \mathbf{\mx})^{\top} \mathbf{J}(\bmthe') (\bfxv - \mathbf{\mx}) + \ln{Z(\bmthe)}.
\end{equation}
Taking into account that the precision matrix $\mathbf{J}(\bmthe')$ is the inverse covariance, the partition function for the Gaussian joint pdf is given by
\[
Z(\bmthe)= (2\pi)^{N/2}\, |\mathrm{det}\mathbf{J}(\bmthe')|^{-1/2} =  (2\pi \lambda)^{N/2} \left[ \mathrm{det} \mathbf{\tilde{J}}(\bmthe'')\right]^{-1/2},
\]
where $\mathbf{\tilde{J}}(\bmthe'')= \la \mathbf{J}(\bmthe')$ is independent of $\la$ [see the definition~\eqref{eq:prec-mat}].  Thus, the \sli NLL is given by
\beq
\label{eq:sli-nll}
-\ln{\mathcal{L}}(\bmthe;\bfxv ) = \frac{1}{2}\left[  (\bfxv - \mathbf{\mx})^{\top} \mathbf{J}(\bmthe') (\bfxv - \mathbf{\mx}) + N\, \ln\la - \ln\mathrm{det}\mathbf{\tilde{J}}(\bmthe'') \right].
\eeq
The above form does not include the constant factor $N\ln(2\pi)/2$ which is irrelevant for the NLL minimization. The trend vector $\mathbf{\mx}$ depends on the parameters $\{ b_{k} \}_{k=1}^{K}$.
The NLL~\eqref{eq:sli-nll} is minimized numerically
using the {\sc Matlab} constrained optimization
function \verb+fmincon+. Constraints (lower and upper bounds) are used to ensure that the parameters are positive and take reasonable values. The  log-determinant in~\eqref{eq:sli-nll}
is calculated numerically using the LU decomposition of the sparse precision matrix.
The optimization parameters include a maximum of $10^4$ iterations and function evaluations, and a tolerance equal to $10^{-4}$ for the cost function and for changes in the optimization variables. The optimization employs the default interior-point method and terminates at  a local minimum after 28 iterations.


\section{Application to Data}
\label{sec:data}
We investigate  \sli-based interpolation for synthetic (simulated) data, reanalysis  \stt data (temperature in degrees Celsius), and ozone measurements over France. The data are   used to provide proof of concept for the \stt-\sli method. In the case of the  synthetic data, we also compare the \sli prediction performance with that of spatio-temporal Ordinary Kriging.

\subsection{Synthetic data}
\label{ssec:synthetic}
\begin{figure}[ht]
  \centering
  \begin{subfigure}[b]{0.6\textwidth}
    \centering\includegraphics[width=\linewidth]{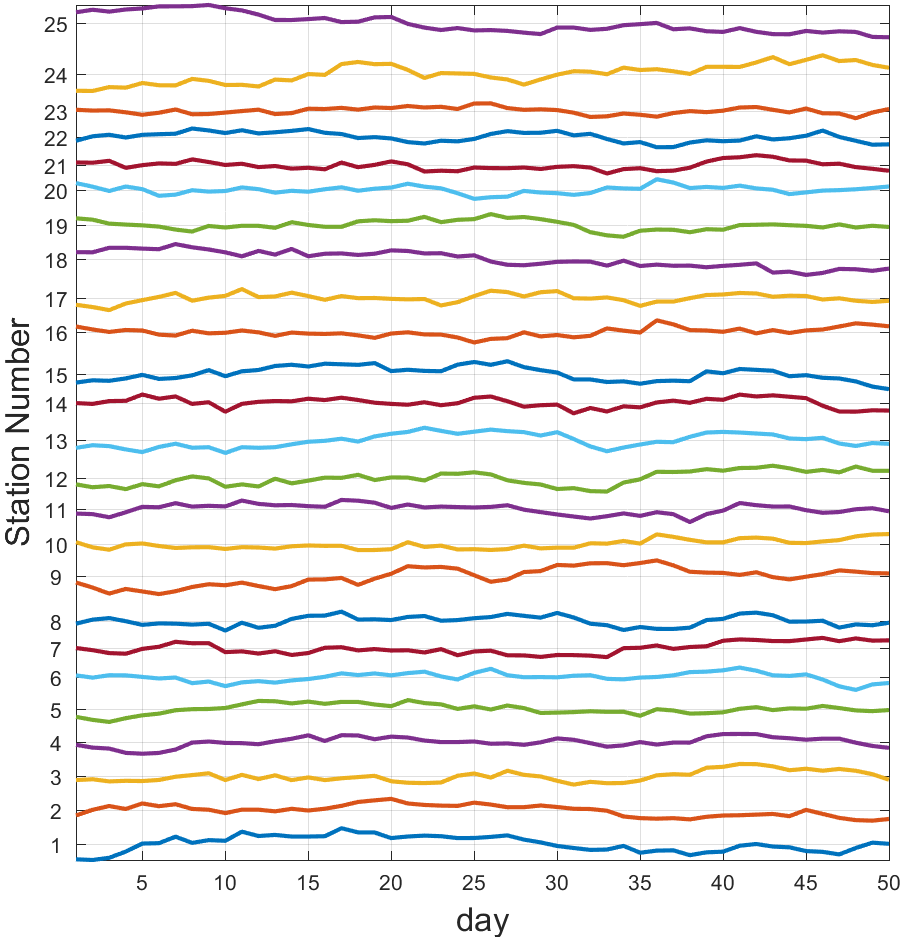}
    \caption{Time series\label{fig:rf4-25days}}
  \end{subfigure} \hfill
  \begin{subfigure}[b]{0.6\textwidth}
    \centering\includegraphics[width=\linewidth]{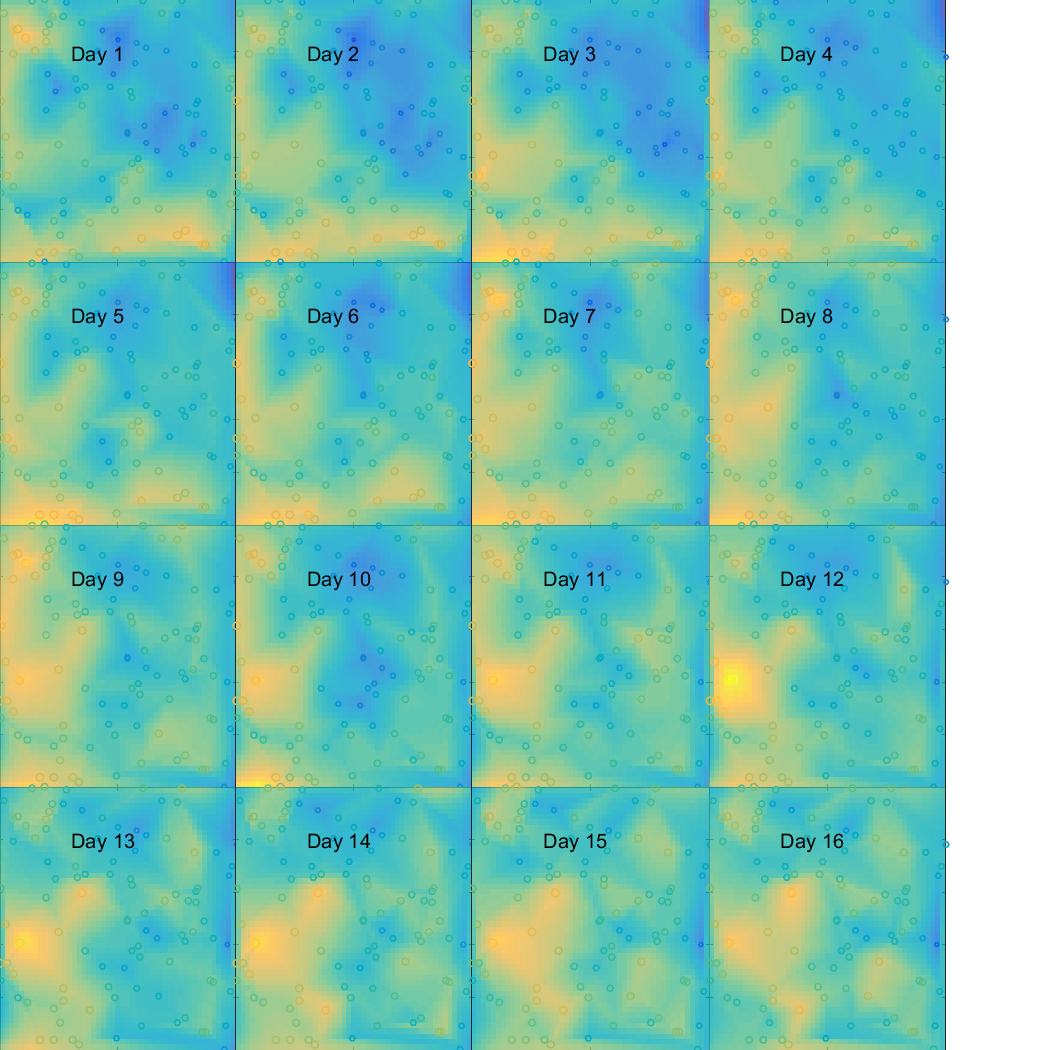}
    \caption{Time slices\label{fig:rf4-slices}}
  \end{subfigure}
  \caption{(a): Time series of synthetic data  at 25 locations. (b): Spatial linear-interpolation maps for the first 16 time slices. The open circles represent the sampling locations. Each time slice corresponds to one simulated day.}
\end{figure}

We generate an \stt realization from a stationary  random field $X(\st;\om)$ with mean $\mx=10$ and variance $\sigma^{2}=5$  using the R package ``RandomFields''~\cite{rf,Schlather15}. The  random field has a separable exponential covariance model with correlation lengths  $\xi{_s}=20$ and $\xi_{t}=10$ in space and time respectively, i.e., $c(\bfr, \tau)= \sigma^{2} \exp(-\rr/\xi_{s}- |\tau|/\xi_{t})$.  The realization is sampled at $N_{s}=100$ random locations over a square spatial domain of length 100 per side, at $N_{t}$ times, i.e., $t_{n}= n\delta t$, for $n=1, \ldots, 50$ where $\delta t=1$. Thus, the sample involves a total of $5000$ points. The resulting time series at 25 spatial locations are shown in Fig.~\ref{fig:rf4-25days}, while the spatial configurations for the first 16 time slices are shown in Fig.~\ref{fig:rf4-slices}.

\subsubsection{\sli parameter estimation using MLE}
We assume that the trend model is simply a constant term, i.e., $b_{1}$. The optimal model parameters are estimated by minimizing the NLL given by~\eqref{eq:sli-nll}.
The orders of the spatial and temporal neighbors  are set to $K_{s}=K_{t}=3$.
The initial guesses for the \sli parameters and the parameter bounds are given in Table~\ref{tab:RF4-slices}.
The value of the cost function (NLL) for the optimal \sli parameters is $\approx -1.4578\times 10^{4}$. The values of the optimal \sli parameters are listed in Table~\ref{tab:RF4-slices}.

\begin{figure}
  \centering
  \begin{overpic}[scale=0.75]{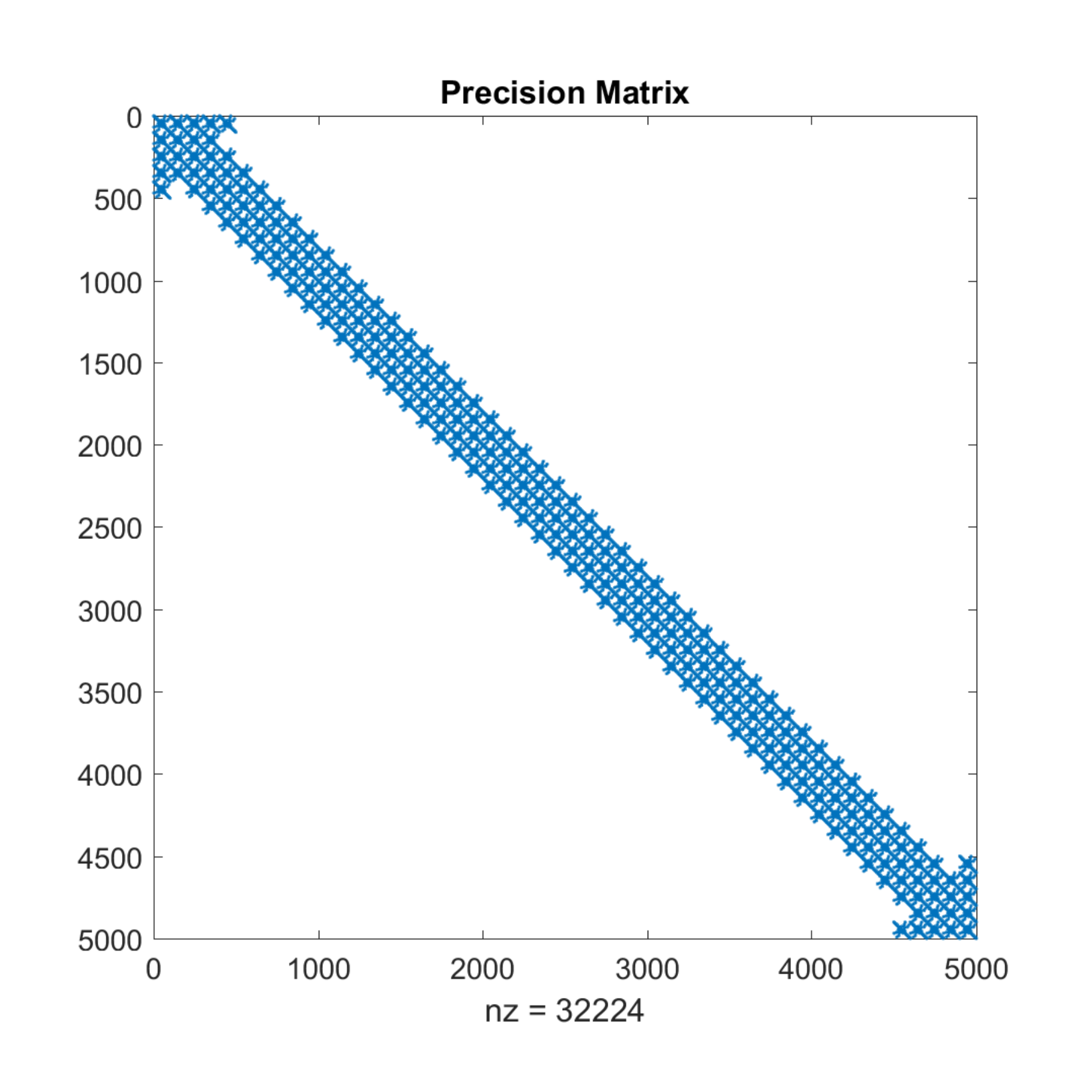}
     \put(15,15){\includegraphics[scale=0.32]{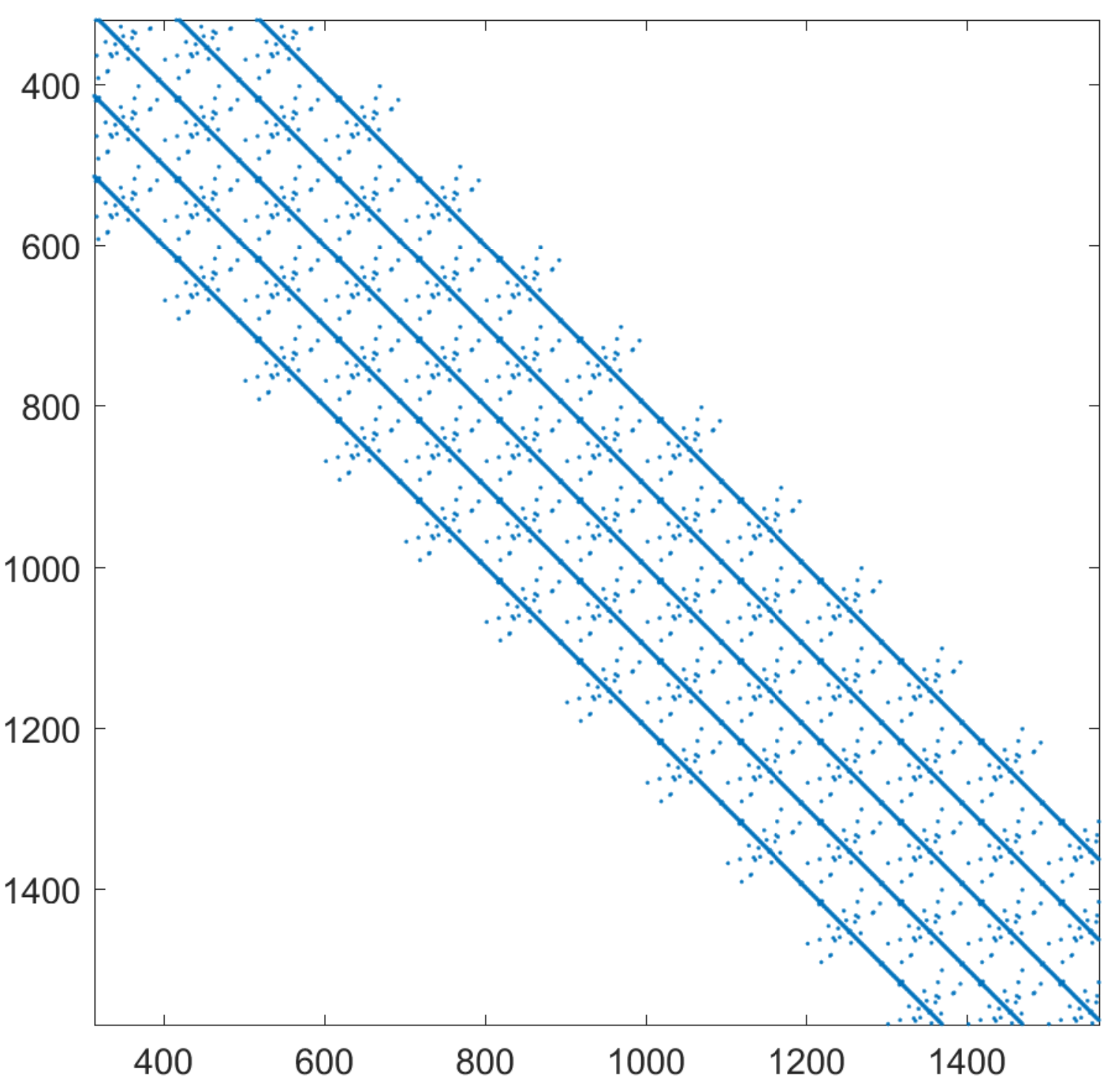}}
  \end{overpic}
 \caption{Sparsity pattern of the \sli precision matrix for the synthetic \stt data; $nz=32\,224$ is the number of non-zero elements in the matrix, leading to a sparsity index of $\approx 0.13\%$. The inset figure shows a detail of the precision matrix. }
 \label{fig:rf4-logabsprec}
\end{figure}

The sparsity pattern of the precision matrix  evaluated with the optimal \sli parameters is shown in Fig.~\ref{fig:rf4-logabsprec}. The non-zero matrix entries are shown as blue dots. The four diagonal bands (two above and two below the main diagonal) comprise nearest and next-nearest temporal neighbors. The sparsity of the precision matrix is evident in the plot; the sparsity index is $\approx 0.13 \%$ corresponding to 32\,224 non-zero elements.


\begin{table}
\centering
\renewcommand*{\arraystretch}{1.2}
\begin{tabular}{lccccc}
\hline
& $\mx$ & $\lambda$ & $c_1$ & $\mu_t$ & $\mu_s$  \\
\hline
Initial values & 9.7424 & 430 & $2.8484\times 10^6$ & 0.5 & 1  \\
\hline
Lower bounds & 9.5798 & $10^{-3}$ & 1 & 1.4 & 0.4 \\
\hline
Upper bounds & 9.9050  & $10^{7}$ & $10^{7}$ & 10 & 10 \\
\hline
Based on MLE & 9.7424 & 0.001  & 4951.7 & 1.4  & 0.40
\\
\hline
\end{tabular}
\caption{\sli parameters for the synthetic \stt data based on MLE. The initial guesses for the optimization are determined by running leave-one-out cross validation using root mean square error as the cost function. The lower and upper bounds on the mean are based on $\overline{x} \mp {5s_{x}}/{\sqrt{N}}$, where $\overline{x}$ is the sample mean,  $s_x$ is the sample standard deviation, and $N=5000$ is the total number of points.}
\label{tab:RF4-slices}
\end{table}

\subsubsection{\sli model performance}
\label{ssec:sli-model-performance}
\begin{figure}[ht]
\centering
\begin{subfigure}[b]{0.6\textwidth}
\centering\includegraphics[width=\linewidth]{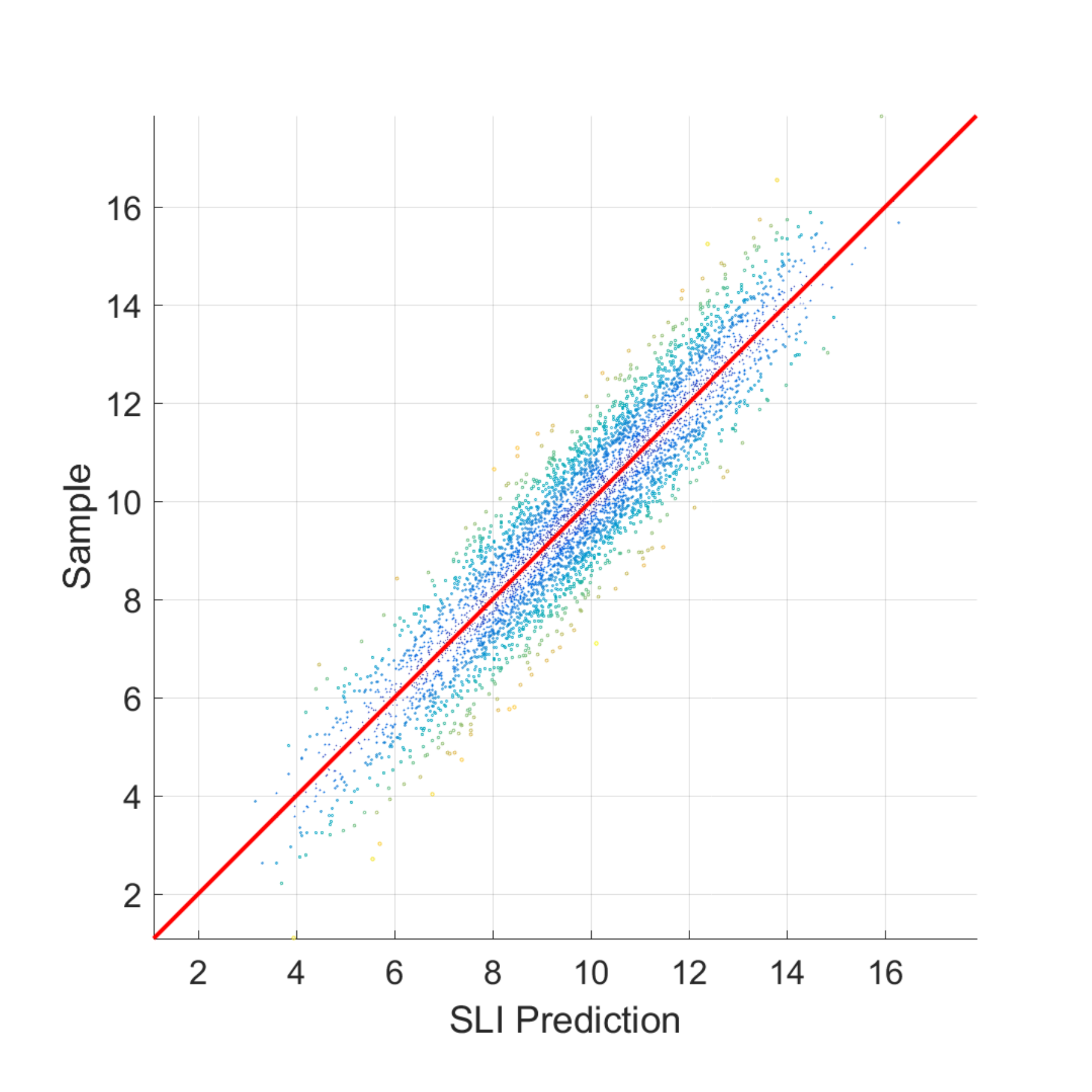}
\caption{Scatter plot\label{fig:rf4-scatter}}
\end{subfigure}  \hfill
\begin{subfigure}[b]{0.6\textwidth}
\centering\includegraphics[width=\linewidth]{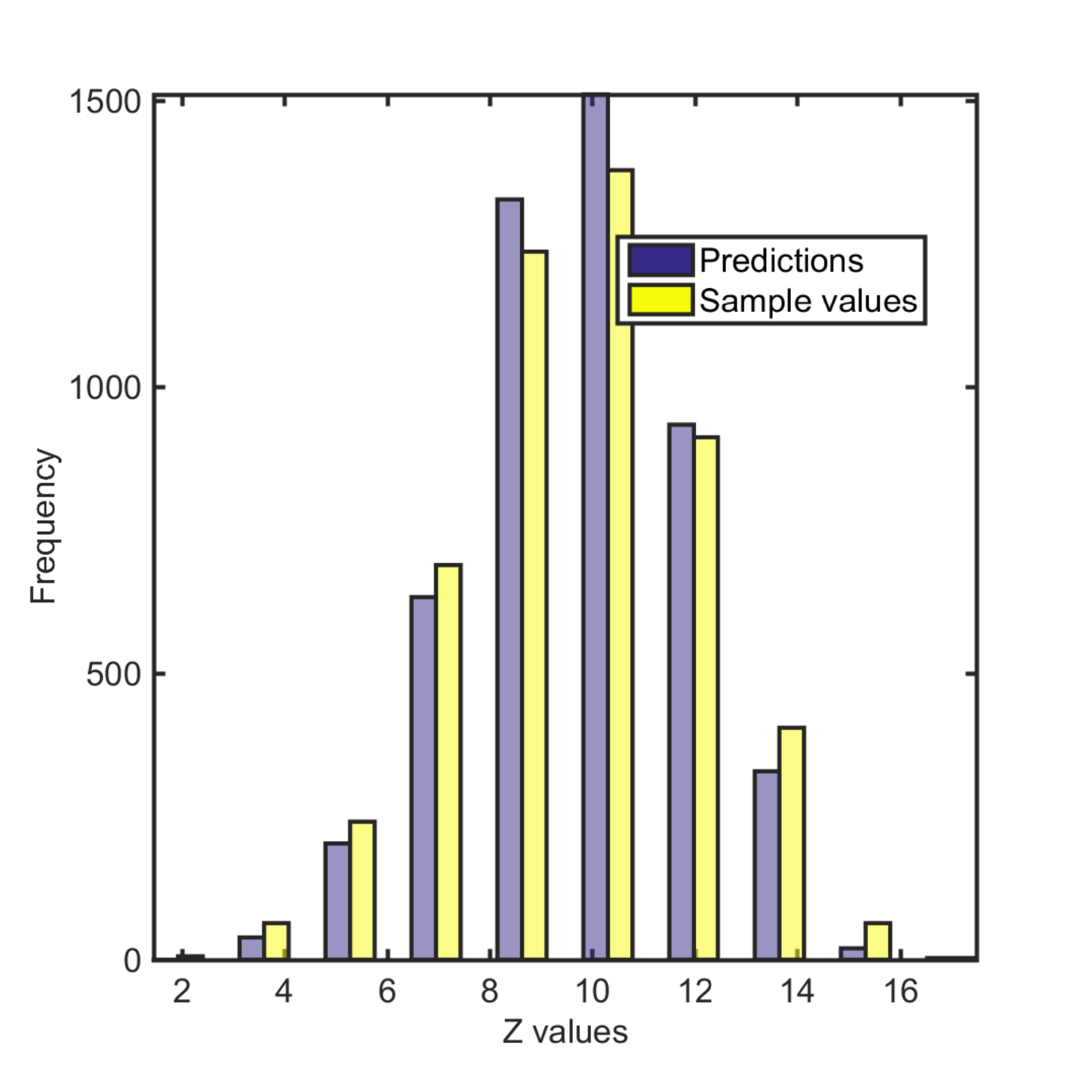}
\caption{Histogram\label{fig:rf4-histo}}
\end{subfigure}
\caption{(a): Scatter plot of the predictions versus the sample values for the synthetic space-time data. (b): Histograms of the sample (yellow) and \sli-predicted (grey) values.}
\end{figure}

To test the performance of the estimated \sli model we use one-slice-out cross validation: we remove and subsequently predict all the values for one time slice using the  sample values at the $N_{t}-1$ remaining time slices. We repeat this experiment by removing sequentially all the time slices, one at a time.
The scatter plot of the predictions (for all $N$ points) versus the sample values is shown in Fig.~\ref{fig:rf4-scatter} and exhibits overall good agreement between the two sets. The histogram plots of the predicted versus the sample values, shown in Fig.~\ref{fig:rf4-histo}, demonstrate that the \sli predictions tend to cluster around the center of the distribution more  than the sample values.

Cross validation measures are presented in Table~\ref{tab:RF4-validation}: ME stands for the mean error (bias), MAE is the mean absolute error, MARE is the mean absolute relative error ($\approx 7.1\%$), RMSE is the root mean square error, RMSRE is the root mean square relative error ($\approx 10.49\%$), $R$ is the linear correlation coefficient ($\approx 0.94$) and $R_{S}$ the Spearman correlation coefficient (also $\approx 0.93$). The validation measures indicate overall good performance of the \sli model with small bias $\approx -0.007$ and very good correlation  $\approx 0.94$. The RMSE is $\approx 0.80$.

In  Table~\ref{tab:RF4-validation} we also compare the SLI cross validation measures with those obtained by means of ST Ordinary Kriging (OK)~\cite{Pebesma16,Wikle19}. The latter is implemented using the function \verb+krigeST+ from the R package  \verb+gstat+. The covariance parameters are estimated by means of the method of moments (MoM), e.g.~\cite{Chiles12}. We report OK cross validation results with three  different parameter sets: The first set (OK-Ex) comprises the parameters of the theoretical covariance function. The second set (OK-Est-1) is based on the optimal covariance model which is fitted to the MoM estimator using unconstrained optimization. Finally, the third set (OK-Est-2) is obtained by means of the same fitting procedure by means of constrained optimization which  forces the model parameters to lie within specific intervals (see caption of Table~\ref{tab:RF4-validation}). The \sli prediction performs better than  OK-Est-1, but it is inferior to OK-Ex and OK-Est-2, while OK-Est-2 has the best performance.  These results are not surprising, given that OK employs the exponential covariance model that was used to generate the data. Nonetheless,  the \sli performance is competitive with that of OK.

\begin{table}
\centering
\renewcommand*{\arraystretch}{1.2}
\begin{tabular}{lrcccccc}
\hline
Method & ME & MAE & MARE & RMSE & RMSRE & $R$ & $R_{S}$ \\
 \hline
SLI & $-$0.0070  &   0.6361  & 0.0717  &  0.7980  &  0.1049   &  0.9383  &  0.9347  \\
 OK-Ex  &  0.0003 &	0.6057  &	0.0684	&   0.7591  &	0.1015 &   0.9444 &   0.9401 \\
 OK-Est-1  &   0.0031  &   0.7807  &   0.0887  &   0.9819  &   0.1326 &   0.9069  &   0.8996  \\
 OK-Est-2 &  $-$0.0011 &   0.5920 &	0.0661 &	0.7398  &	0.0918 &	0.9468  &	0.9425 \\
\hline
\end{tabular}
\caption{One-slice-out cross validation (CV) interpolation performance for the Gaussian data with separable exponential covariance. The CV measures are calculated by comparing  the true values of each time slice (from 1 to 50) and the predicted values. The  predictions are based on $N_{t}-1$ time slices that exclude  the predicted slice. First row: Predictions based on SLI. Second row (OK-Ex): Predictions based on Ordinary Kriging with theoretical covariance parameters.  Third row (OK-Est-1): Predictions based on Ordinary Kriging with estimated covariance parameters (unconstrained estimates). Fourth row (OK-Est-2): Predictions based on Ordinary Kriging with estimated covariance parameters using constraints: $\xi_{s} \in [10, 30]$, $\xi_t \in [5, 15]$, nugget variance $\in  [0, 0.1]$, $\sigma^2 \in [2.5, 7.5]$.}
\label{tab:RF4-validation}
\end{table}

\subsection{Hourly temperature reanalysis data}
\label{ssec:era5}
\begin{figure}[ht]
  \centering
  \begin{subfigure}[b]{0.5\textwidth}
    \centering\includegraphics[width=\textwidth]{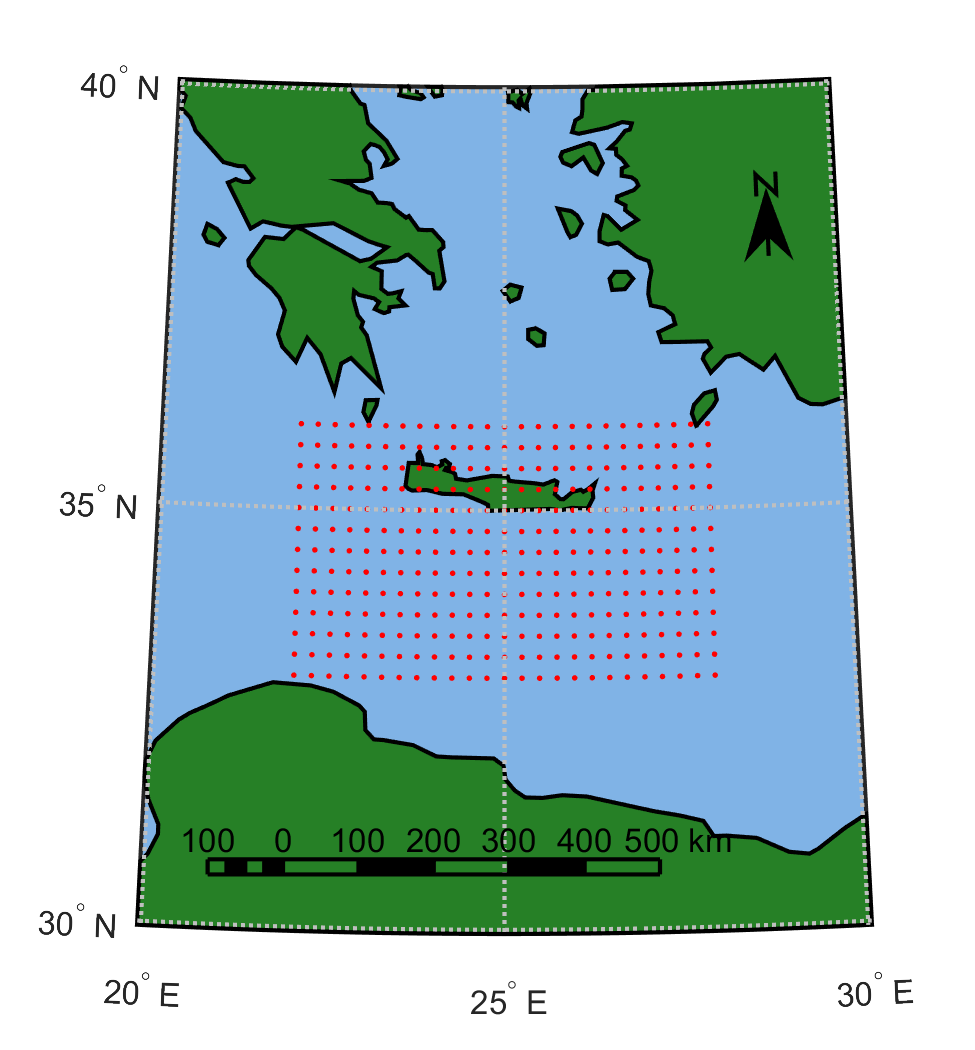}
    \caption{Map showing the grid sites (red markers) of the ERA5 reanalysis data around the island of Crete that are used in the temperature analysis in  Section~\ref{ssec:era5}.  \label{fig:era-grid}}
  \end{subfigure}  \hfill
    \begin{subfigure}[b]{0.6\textwidth}
    \centering\includegraphics[width=\linewidth]{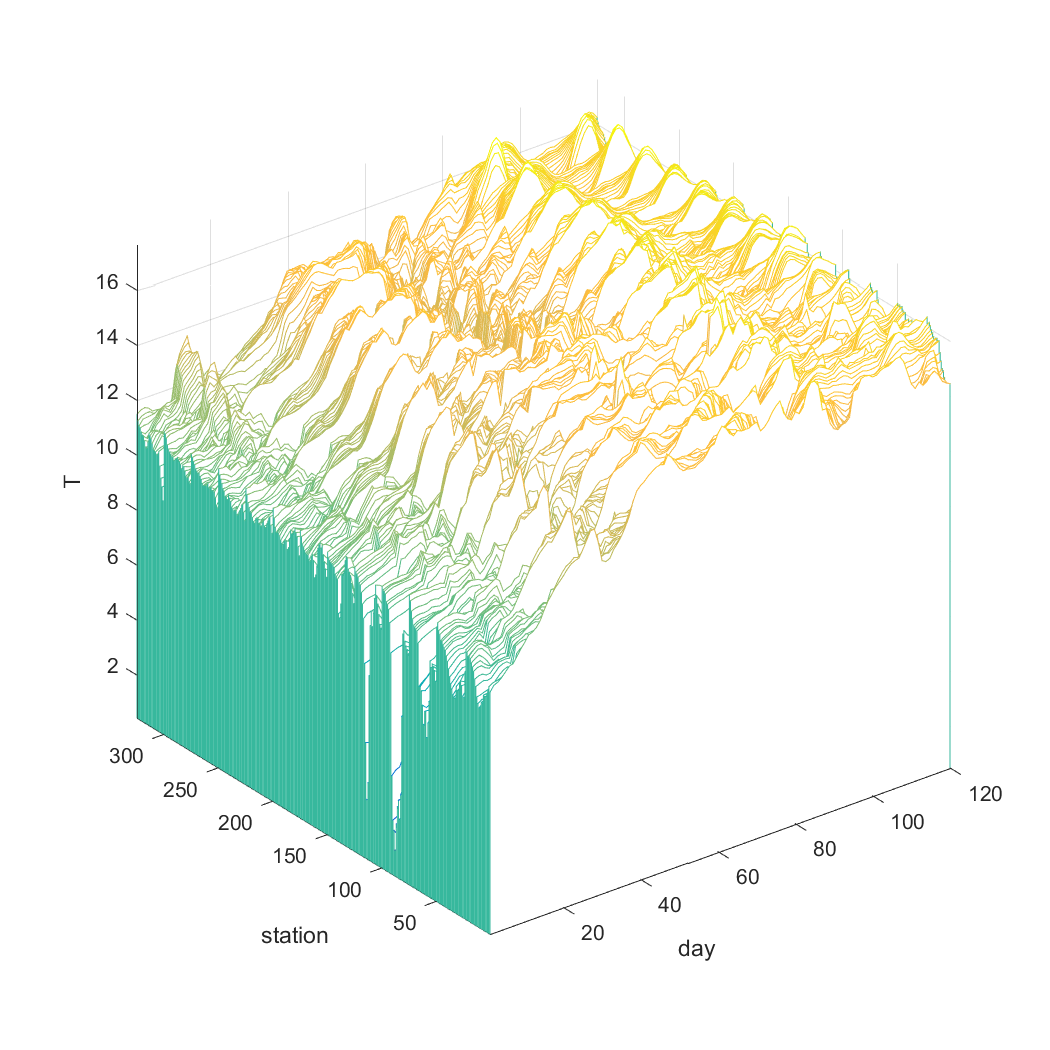}
    \caption{Temperature \label{fig:era-time-series}}
  \end{subfigure}
  \caption{(a): Spatial grid for the ERA5 temperature data (degrees  Celsius) around the island of Crete (Greece). (b): Time series of  temperature (in degrees Celsius) at the ERA5 grid sites shown in (a).}
\end{figure}

We use ERA5 reanalysis temperature data (degrees Celsius) downloaded from the Copernicus Climate Change Service~\cite{copernicus18}.  The dataset  includes 39\,000 points that correspond to hourly  values for five consecutive days (January 1-5, 2017) at the nodes of a
$13 \times 25$ spatial grid around the island of Crete (Greece) as shown in Fig~\ref{fig:era-grid}. The average spatial resolution is $\approx0.28$ degrees (grid cell size $\approx 31$km). The data are displayed as time series in Fig.~\ref{fig:era-time-series}.

\begin{figure}
\centering    \includegraphics[width=.89\textwidth]{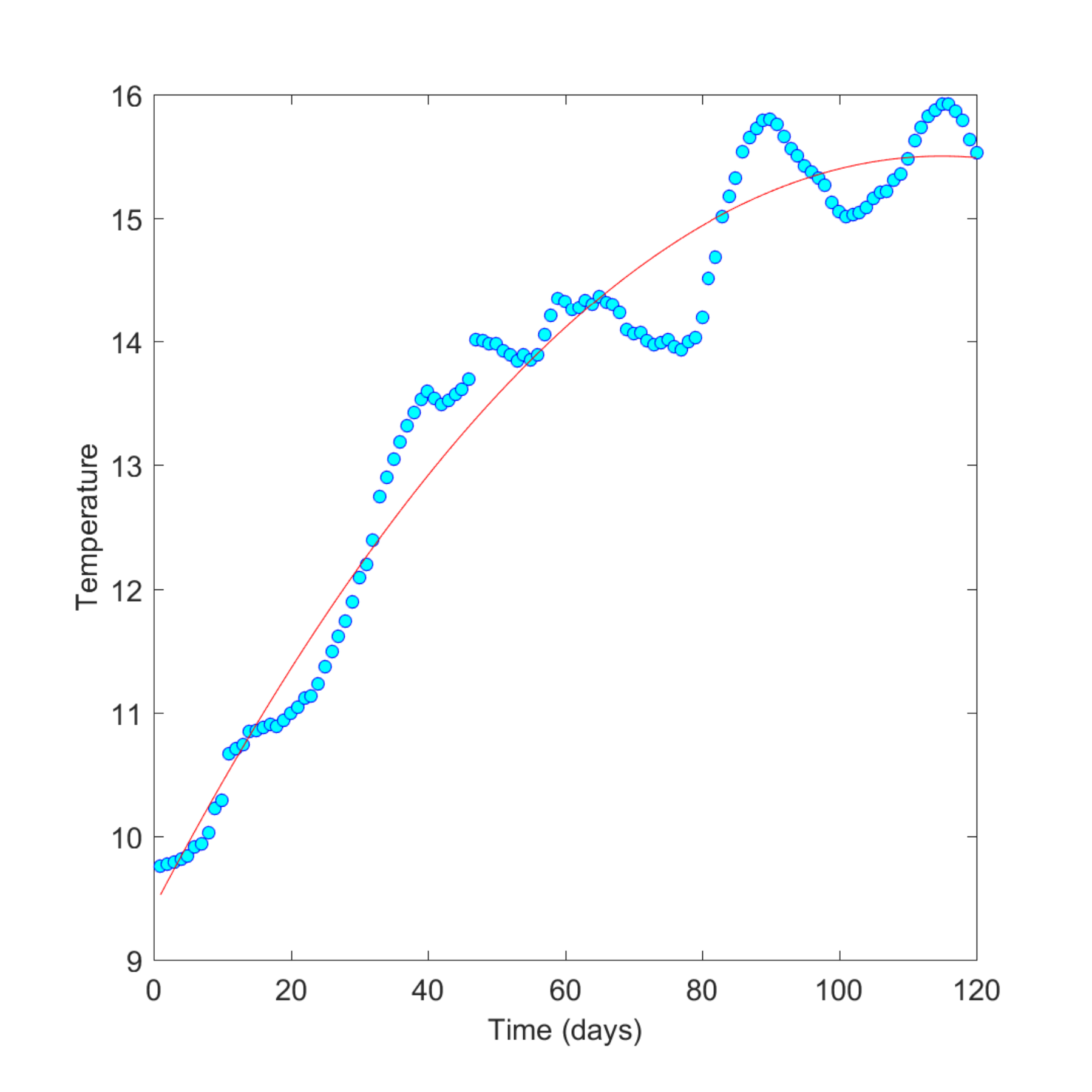}
\caption{Plot of the spatially averaged temperatures (in degrees Celsius) as a function of time and least-squares fit to a second-degree polynomial function.}
\label{fig:era5-trend}
\end{figure}

The temperature data exhibit a clear increasing trend in time during the studied period. This is evidenced in the plot of the spatially averaged temperature as a function of time in Fig.~\ref{fig:era5-trend},  and the temperature fit  with the linear regression model $\mx(t)=b_{1}+b_{2} t + b_{3} t^{2}$ (where $t$ is measured in days). Thus, the parameter vector~\eqref{eq:bmthe} with $K_{s}=K_{t}=3$ is given by
$\bmthe=\left( b_{1}, b_{2}, b_{3}, \la, c_{1}, \mu_{s}, \mu_{t}, 3, 3 \right)^\top$.

The \sli parameter estimation and the performance assessment are carried out as in the  synthetic data case study (Section~\ref{ssec:synthetic}). The parameter estimates are shown in Table~\ref{tab:ERA-slices}.
The  precision matrix has a sparsity index $\approx 0.006\%$, corresponding to  $931\,070$ non-zero  entries out of $1.521\times 10^9$ entries.

\begin{figure}[ht]
  \centering
  \begin{subfigure}[b]{0.6\textwidth}
    \centering\includegraphics[width=\linewidth]{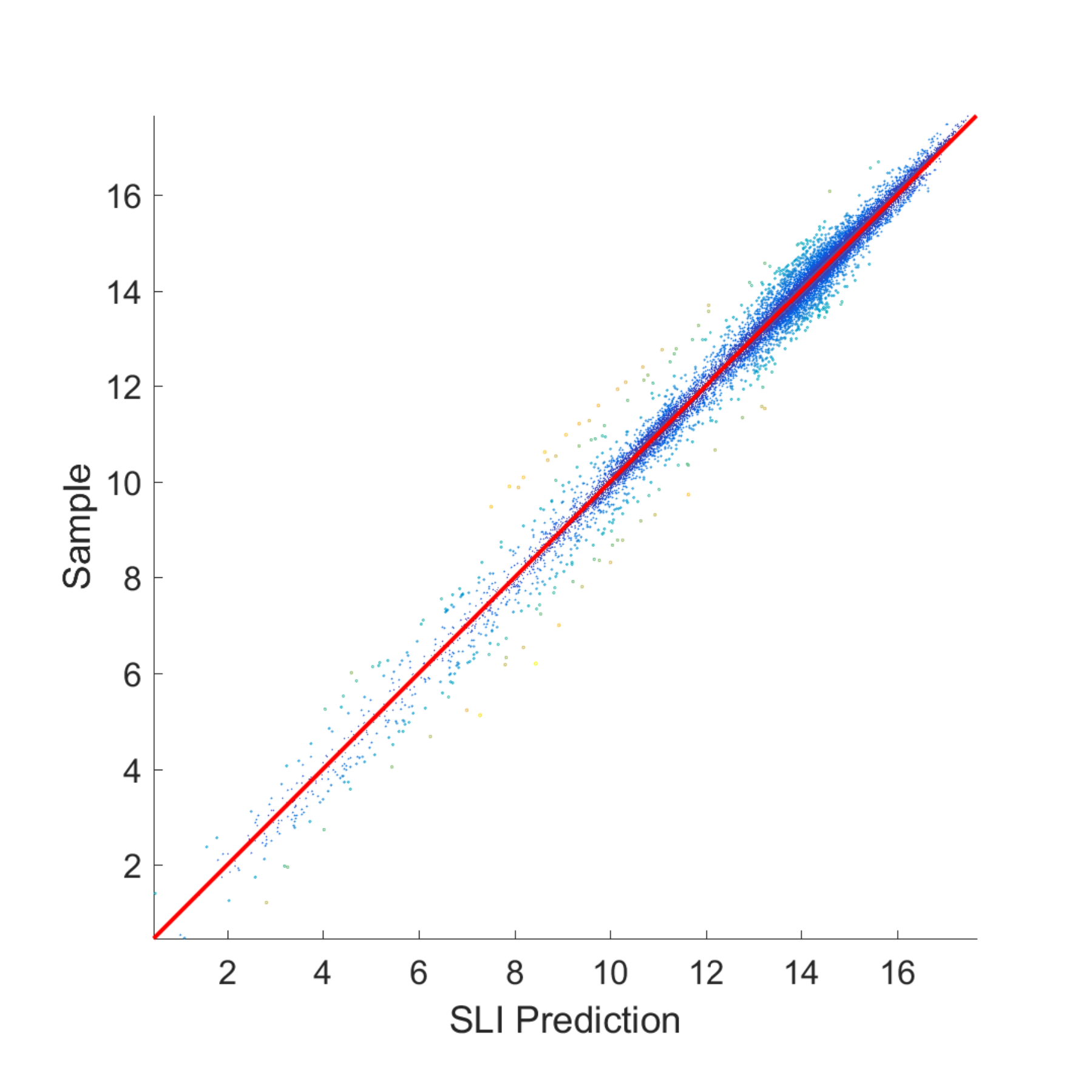}
    \caption{Scatter plot\label{fig:era-scatter}}
  \end{subfigure}  \hfill
    \begin{subfigure}[b]{0.6\textwidth}
    \centering\includegraphics[width=\linewidth]{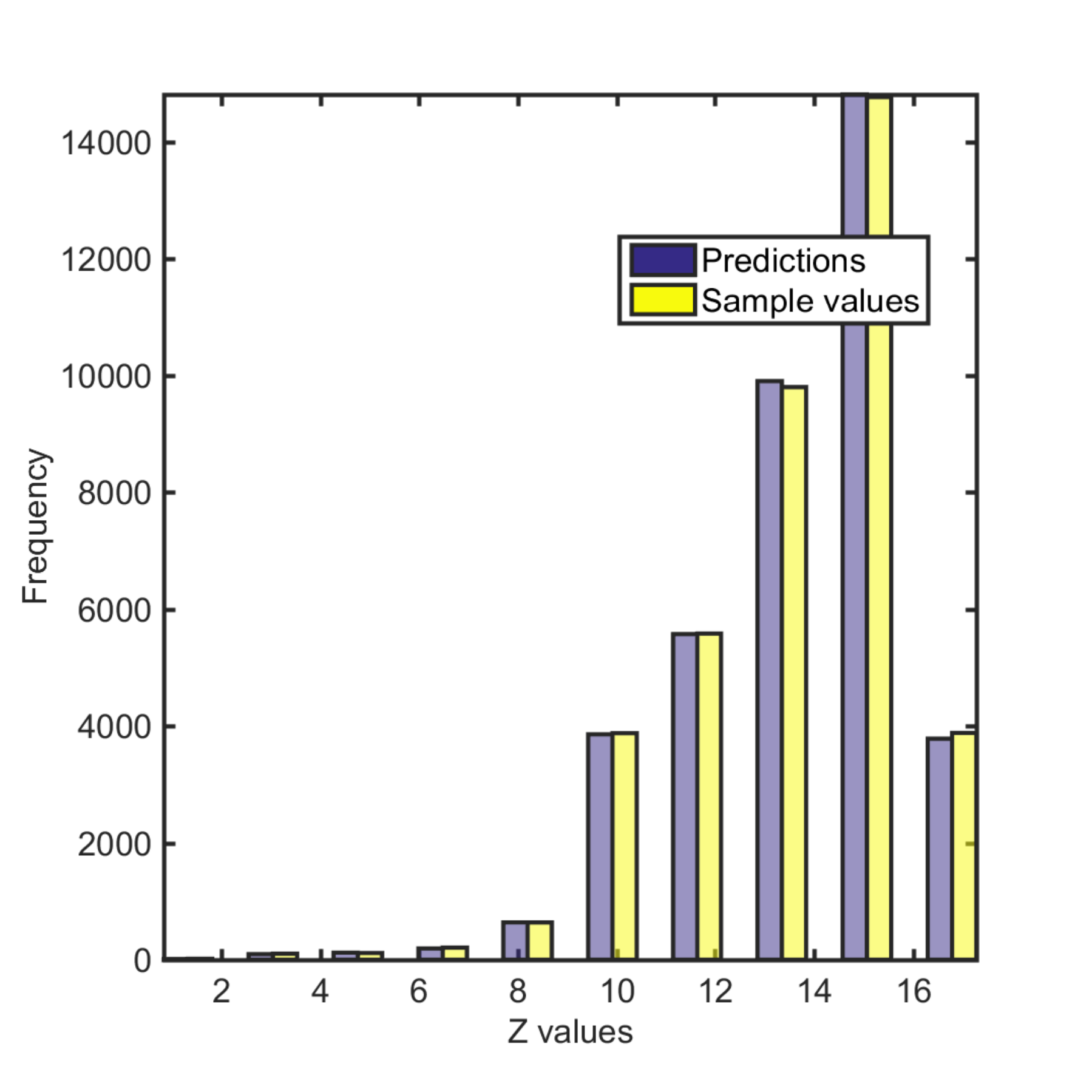}
    \caption{Histogram\label{fig:era-histo}}
  \end{subfigure}
  \caption{(a): Scatter plot of the \sli predictions versus the sample values for the ERA5 temperature data. (b): Histograms of the sample (yellow) and \sli-predicted (grey) values.}
\end{figure}

The scatter plot of the predictions (for all $N$ points)  is shown in Fig.~\ref{fig:era-scatter} and exhibits good agreement between the predictions and the data. The histogram plots of the predicted versus the sample values (Fig.~\ref{fig:era-histo}) also show that \sli predictions have lower dispersion than the sample values, as was the case for the synthetic data. The cross validation measures (obtained by sequentially removing each of the 120 hourly time slices) are shown in Table~\ref{tab:ERA-T-validation} and  confirm the interpolation performance for the \sli model.

\begin{table}
\centering
\renewcommand*{\arraystretch}{1.2}
\begin{tabular}{lccccccc}
\hline
& $b_{1}$ & $b_{2}$ & $b_{3}$ & $c_1$ & $\mu_t$ & $\mu_s$ & $\lambda$ \\
\hline
Initial  & 9.4203 & 0.1058 & $-4.6 \times 10^{-4}$ & 300 & 3  & 2.5 & 10 \\
\hline
L.B. & 9.2042 & 0.0975 & $-5.26 \times 10^{-4}$ & 1 & 1  & 1 & 1 \\
\hline
U.B. & 9.6363 &  0.1140 & $-3.94 \times 10^{-4}$ & 1000 & 1000  & 10 & 10 \\
\hline
MLE & 9.5587 & 0.1019 & $-0.0004$ & 111.4797 & 1 & 1 & 1 $\times 10^{-4}$\\
\hline
\end{tabular}
\caption{\sli model parameters for the ERA5 temperature data based on MLE.
The initial values for $b_{1}, b_{2}, b_{3}$  are obtained from the  coefficients of the regression model for the trend.  The lower (L.B.) and upper (U.B.) bounds of the coefficients are the respective limits of the regression-based 95\% confidence intervals.}
\label{tab:ERA-slices}
\end{table}

\begin{table}
\centering
\renewcommand*{\arraystretch}{1.2}
\begin{tabular}{ccccccc}
\hline
 ME ($^\circ$C) & MAE ($^\circ$C) & MARE & RMSE ($^\circ$C) & RMSRE & $R$ & $R_{S}$ \\
 \hline
 $-0.0008$ &   0.1022  &  0.0085  &  0.1737  &  0.0216  &  0.9971  &  0.9969  \\
\hline
\end{tabular}
\caption{One-slice-out cross validation (CV) test of the \sli interpolation performance for the ERA5 temperature data. The CV measures are calculated by comparing  the true temperature values of each hourly time slice (from 1 to 120) and the \sli predictions that are based on the \sli model with the MLE parameters reported in Table~\ref{tab:ERA-slices}. The predictions are based on $N_{t}-1$ time slices excluding  the predicted slice.}
\label{tab:ERA-T-validation}
\end{table}

\subsection{Hourly ozone concentration data}
\begin{figure}[ht]
  \centering
  \begin{subfigure}[b]{0.6\textwidth}
    \centering\includegraphics[width=\linewidth]{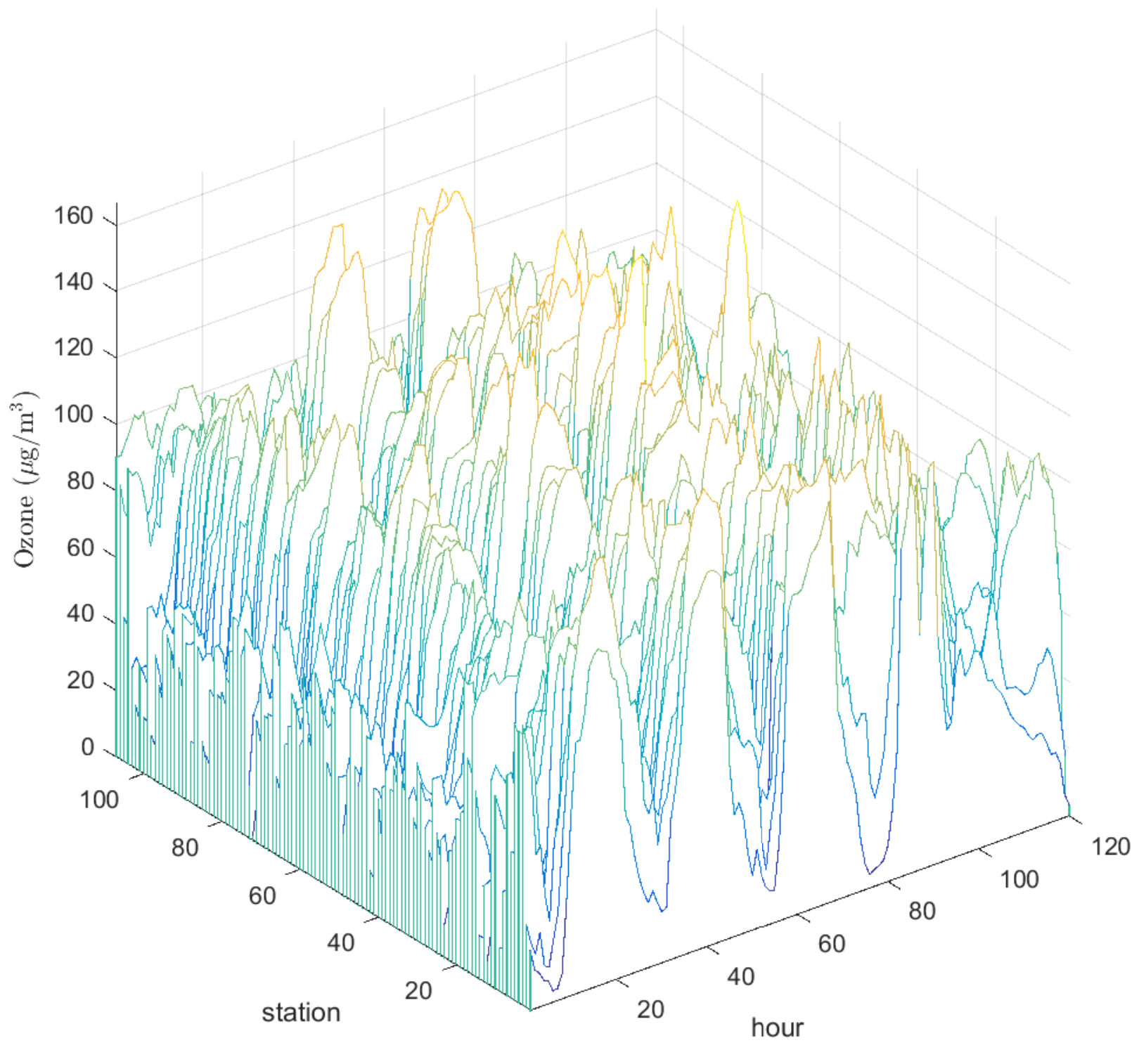}
    \caption{Time series\label{fig:TS_rfsummer}}
  \end{subfigure} \hfill
  \begin{subfigure}[b]{0.6\textwidth}
    \centering\includegraphics[width=\linewidth]{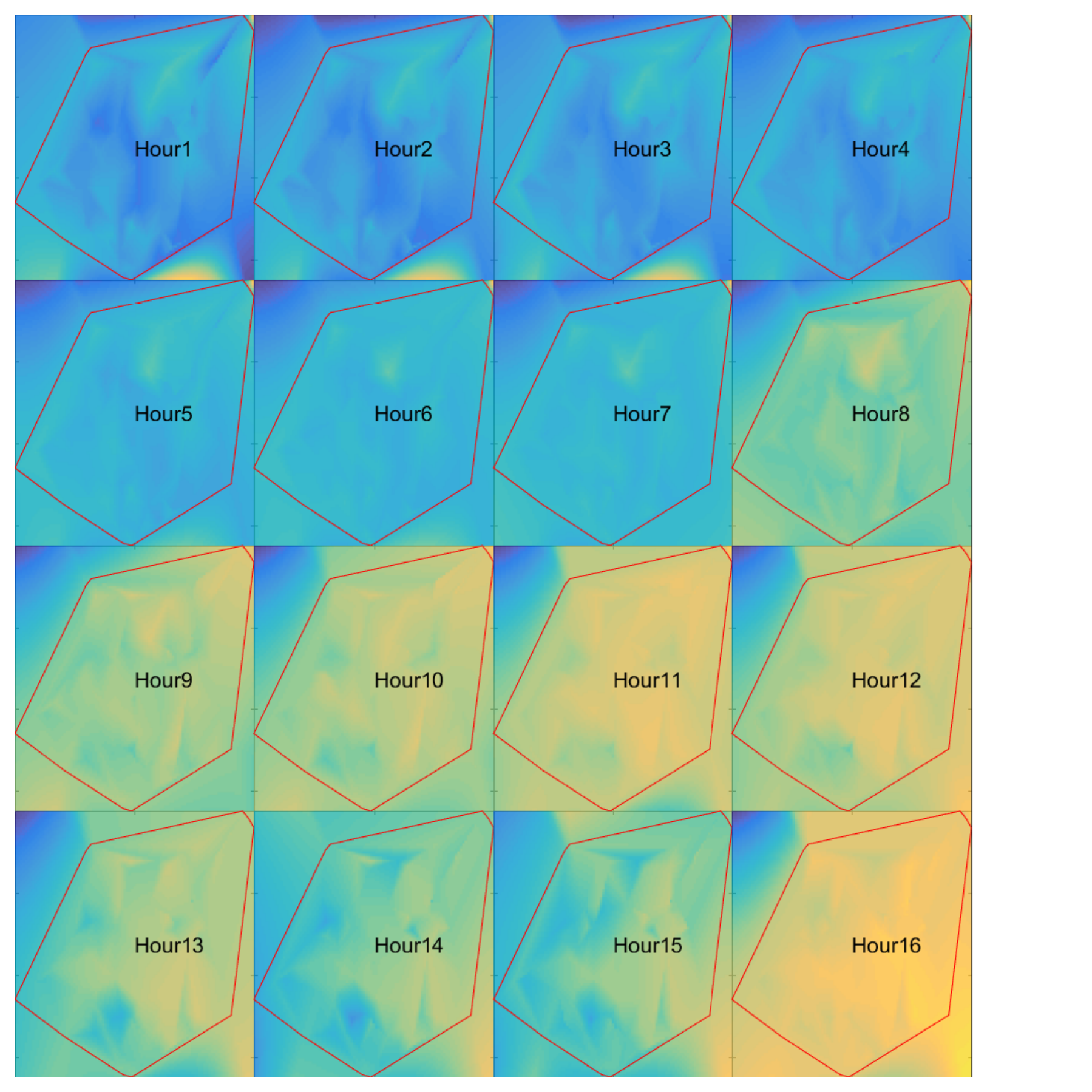}
    \caption{Time slices\label{fig:rfsummer-slices}}
  \end{subfigure}
  \caption{(a): Time series of ozone O\textsubscript{3} hourly data ($\mu$g/m\textsuperscript{3}) at 107  stations in France and spatial linear-interpolation maps for the first 16 hourly time slices. (b): Linear interpolation maps of ozone concentration  for different time slices (hours of the day). The red polygon marks the convex hull of the station network.}
\end{figure}

This dataset includes ozone (O\textsubscript{3})  hourly concentration data (measured in $\mu$g/m\textsuperscript{3}) for five consecutive days (July 1-5, 2014), downloaded from the French  GEOD'AIR database (web site: \url{www.prevair.org}). The data are collected at 335 scattered stations distributed around France. The time series at 107 stations that reported data at all times are shown in Fig.~\ref{fig:TS_rfsummer}. Linear interpolation maps for the first 16 time slices are shown in Fig.~\ref{fig:rfsummer-slices}.

A visual inspection of the spatially averaged ozone time series indicates the existence of a temporal trend which is modeled by means of the following function which exhibits daily (24-hr) periodicity
\beq
\label{eq:time-trend}
\mx(t)=b_{1}+ \left(b_{2}  + b_{3} t + b_{4} t^{2}\right) \cos\left(\frac{2\pi t}{24} \right) +
\left(b_{5}  + b_{6} t + b_{7} t^{2}\right) \sin\left(\frac{2\pi t}{24} \right).
\eeq

\begin{figure}
\centering
\includegraphics[width=.95\textwidth]{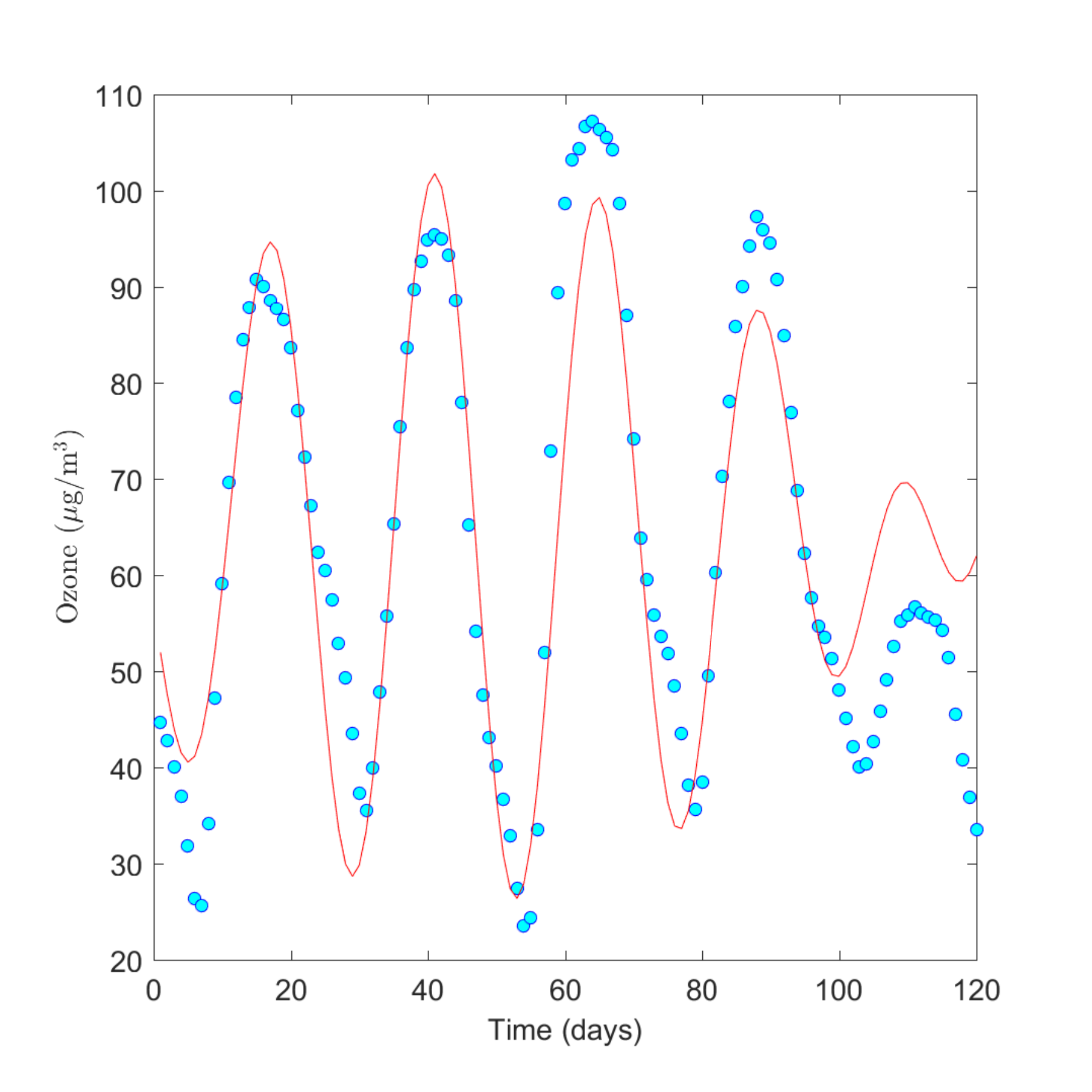}
\caption{Plot of the spatially averaged ozone concentration ($\mu$g/m\textsuperscript{3}) over the study domain as a function of time (circles) and least-squares fit to the temporal trend  function given by~\eqref{eq:time-trend}.}
\label{fig:ozone-trend}
\end{figure}

The \sli parameter estimation is conducted using MLE. The parameter estimates are shown in Table~\ref{tab:Horaire_summer-slices}.
MARE and RMSRE are infinite because the dataset includes zero values.
The  precision matrix has a sparsity index $\approx 0.07\%$, i.e., it includes $1.154 \times 10^{5}$ non-zero  entries out of $1.649\times 10^8$ entries.

\begin{table}
\centering
\renewcommand\tabcolsep{3pt} 
\renewcommand\arraystretch{1.2} 
\resizebox{\textwidth}{!}{
\begin{tabular}{lccccccccccc}
\hline
& $b_{1}$ & $b_{2}$ & $b_{3}$ & $b_{4}$ & $b_{5}$ & $b_{6}$ & $b_{7}$ & $c_1$ & $\mu_t$ & $\mu_s$ & $\lambda$ \\
\hline
Initial  & 64.05 & $-7.19$ & $-0.14$ & $1.5\times10^{-3}$ &
$-18.59$ & $-0.77$ & $8\times 10^{-3}$ & 1 & 3  & 2.5 & 1 \\
\hline
L.B. & 62.13 & $-15.50$ & $-0.45$ & $-10^{-3}$ &
$-26.84$ & $-1.08$ & $5.6\times10^{-3}$ & $10^{-3}$ & 0.5  & 0.5 & $10^{-3}$ \\
\hline
U.B. & 65.97 & $1.12$ & $0.17$ & $4\times 10^{-3}$ &
$-10.34$ & $-0.45$ & $0.011$ & $10^7$ & 10  & 10  & $10^{7}$ \\
\hline
MLE  & 63.64 & $-7.67$ & $-0.14$ & $1.6\times10^{-3}$ &
$-19.57$ & $-0.57$ & $5.6\times10^{-3}$ & 77.63 & 1.17  & 0.5 & $10^{-3}$\\
\hline
\end{tabular}}
\caption{\sli model parameters for the French ozone concentration data based on MLE.
The initial values for the trend coefficients $\{ b_{i}\}_{i=1}^{7}$  are obtained from the  coefficients of the regression model for the trend~\eqref{eq:time-trend}.  The lower (L.B.) and upper (U.B.) bounds of the coefficients are the limits of the respective regression-based 95\% confidence intervals.  The value of the cost functional (NLL) at the optimum is equal to $-5.471\times 10^{-3}$. }
\label{tab:Horaire_summer-slices}
\end{table}

\begin{figure}[ht]
  \centering
  \begin{subfigure}[b]{0.6\textwidth}
    \centering\includegraphics[width=\linewidth]{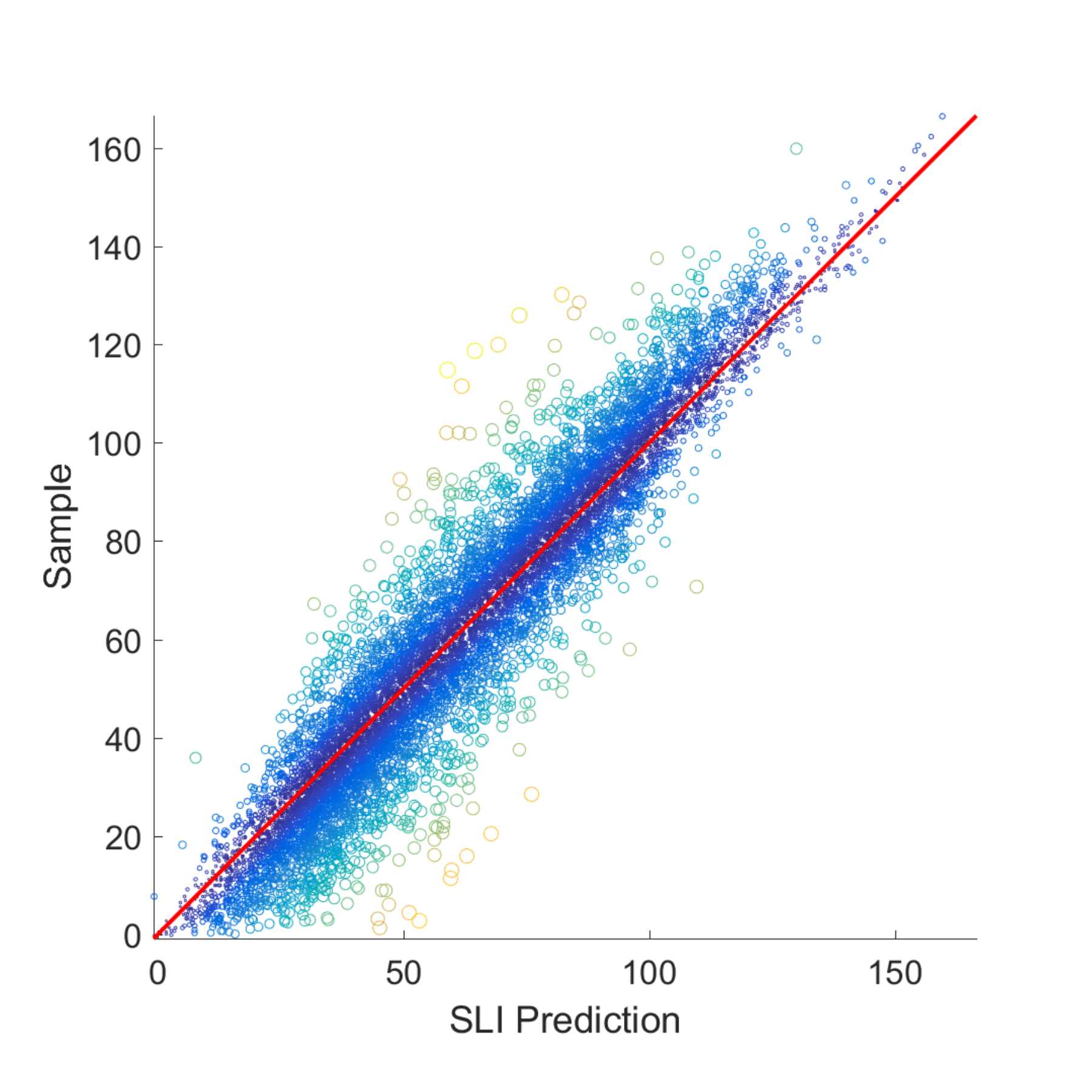}
    \caption{Scatter plot\label{fig:Horaire_summer-scatter}}
  \end{subfigure}  \hfill
    \begin{subfigure}[b]{0.6\textwidth}
    \centering\includegraphics[width=\linewidth]{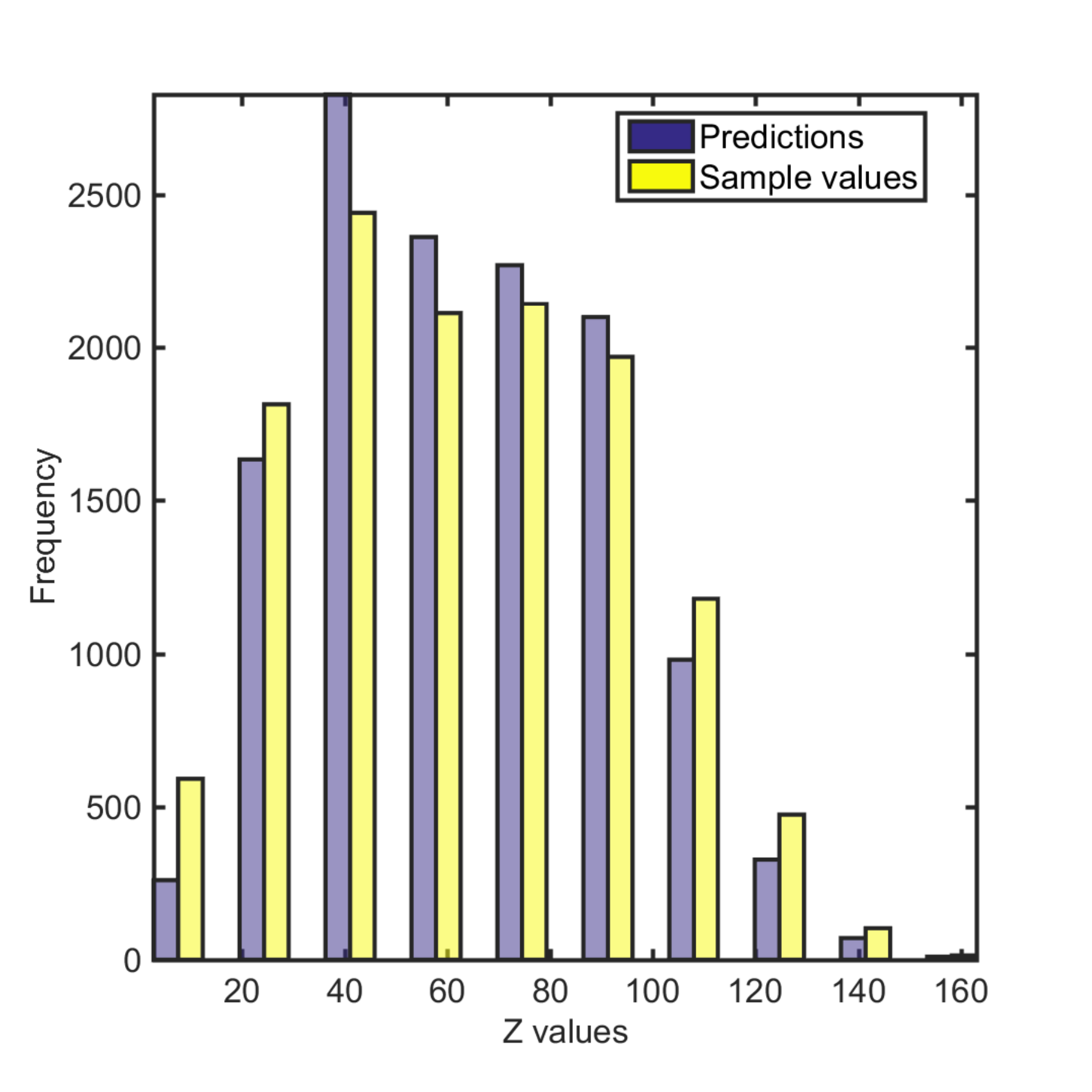}
    \caption{Histogram\label{fig:Horaine_summer-histo}}
  \end{subfigure}
  \caption{(a): Scatter plot of the \sli predictions versus the sample values for the French hourly ozone data. (b): Histograms of the sample (yellow) and \sli-predicted (grey) values.}
\end{figure}

The scatter plot of the predictions (for all $N$ points) versus the sample values is shown in Fig.~\ref{fig:Horaire_summer-scatter} and exhibits overall good agreement between the data and the predictions. The histogram plots of the predicted versus the sample values, shown in Fig.~\ref{fig:Horaine_summer-histo} also show that \sli predictions have lower dispersion than the sample values as in the case of synthetic data. The cross validation performance measures (obtained by sequentially removing each of the 120 hourly time slices) are shown in Table~\ref{tab:Horaire_summer-T-validation}, and they demonstrate very good interpolation performance for the \sli model. Note that $\mu_s = 0.5$, which implies that the spatial bandwidth is small. On the other hand,  $\mu_{t}\approx 1.17$ and $K_{t}=3$ imply that the temporal bandwidth is $\mu_{t} (K_{t}-1) \delta t\approx 2.34$~hr ($\delta t =1$ hour). This result implies that the \sli predictions are at most locations based on the four temporal nearest neighbors (two forward and two backward). For the first and last time slices the bandwidth is $K_{t}\, \mu_{t} \delta t= 3.50$~hr.

\begin{table}
\centering
\renewcommand*{\arraystretch}{1.2}
\begin{tabular}{ccccccc}
\hline
 ME ($\mu$g/m\textsuperscript{3}) & MAE ($\mu$g/m\textsuperscript{3}) & MARE & RMSE ($\mu$g/m\textsuperscript{3}) & RMSRE & $R$ & $R_{S}$ \\
 \hline
 0.00853 &   6.0941  &  Inf  &  8.503  &  Inf  &  0.96  &  0.96  \\
\hline
\end{tabular}
\caption{Cross validation (CV) interpolation performance measures for the ozone concentration values. The CV measures  compare the true ozone concentration values ($\mu g/m^{3}$) of each hourly time slice (from 1 to 120) with the \sli predictions. The latter are based on $N_{t}-1$ time slices that exclude  the predicted slice.}
\label{tab:Horaire_summer-T-validation}
\end{table}

\section{Discussion and Conclusions}
\label{sec:conclusions}
We present a theoretical framework for constructing \stt models based on exponential Boltzmann-Gibbs joint probability density functions. The \stt-\sli model presented herein exploits   an energy function with local interactions which imposes sparse structure on the precision matrix.
 The local interactions are implemented by means of compactly supported kernel functions that compensate for the lack of a structured lattice.
 However, the model is also applicable to regular lattice data.
 In this case the \sli model is equivalent to a Gauss Markov random field with a specific precision matrix structure. To our knowledge, this structure that involves kernel-matrix weights has not been used before in models of real-valued \stt data. The \stt-\sli model extends the purely spatial \sli model~\cite{dth15} to the space-time domain. The \sli approach  shares the reliance on kernel functions with kernel-based  reconstruction of graph signals~\cite{Romero17}.  In the case of environmental \stt data, the graph topology is not given \emph{a priori}, but it is determined from the data by optimizing the negative likelihood of the model.

 The \sli model presented features a Gaussian energy functional with a  sparse precision matrix.   Explicit expressions are given for \stt prediction and the conditional variance at the prediction sites. The sparse precision matrix representation allows computationally efficient implementation of parameter estimation and prediction procedures. The computational efficiency stems from the fact that in \sli it is not necessary to store and invert large and dense covariance matrices.

The optimization of the cost function (the negative logarithmic likelihood) was based on the interior-point  algorithm which terminates at local minima. The  landscape of the cost function should be further investigated in order to understand the patterns of local minima.  It is also possible to run a global optimization algorithm to search for the global optimum of the cost function. On the other hand, experience with the purely spatial  \sli model~\cite{dth15} shows that local minima of the cost function provide parameter estimates that are sufficient for interpolation purposes.

In terms of prediction performance, we have shown (see Table~\ref{tab:RF4-validation}) that for synthetic data the \sli cross validation statistics are slightly inferior but  competitive with those obtained by means of Ordinary kriging. This behavior is observed in spite of the fact that Ordinary kriging has the advantage of employing the functional form (exponential) of the covariance model used to generate the data. In our experience with performance comparisons based on spatial data, the ranking of different methods with regard to prediction performance may change depending on the specific data set. In our opinion, the results shown herein  establish that \sli is a competitive method for space-time data interpolation. Further studies can elaborate on the performance of \sli relative to other methods.

The formulation presented herein can be extended to multivariate random fields by suitable selection of the energy function.
In addition, it is possible to include  anisotropic and more general (e.g., geodesic) spatial distance metrics in the kernel functions, periodic patterns (in space and in time) by adding shifted averaged squared increments, and spatial dependence of the coefficients $\la$ and $c_{1}$. Such extensions will enhance the flexibility of the \sli model at the cost of some loss in computational efficiency.



\section*{Acknowledgments}
This work was funded by the Operational Program ``Competitiveness, Entrepreneurship and Innovation 2014-2020'' (co-funded by the European Regional Development Fund) and managed by the General Secretariat of Research and  Technology, Ministry of Education, Research and Religious Affairs under the project DES2iRES (T3EPA-00017) of the ERAnet,  ERANETMED\_NEXUS-14-049. This support is gratefully acknowledged.

We thank  Prof. Valerie Monbet (Universit\'{e} de Rennes) for  suggesting the ERA5 reanalysis data and Dr. Denis Allard (INRA)  for the French ozone data collected by the Laboratoire Central de Surveillance de la Qualit\'{e} de l'Air.
Dr. Emmanouil Varouchakis (Technical University of Crete)  helped with data analysis in R.  Prof. Ioannis Emiris (University of Athens) made useful suggestions regarding the computation of the log-determinant of the precision matrix. Finally, we acknowledge  two anonymous reviewers whose comments helped to improve this manuscript overall.

\clearpage

\bibliographystyle{unsrt}

\begin{thebibliography}{10}

\bibitem{National13}
National~Research Council et~al.
\newblock {\em Frontiers in Massive Data Analysis}.
\newblock National Academies Press, Washington, DC, 2013.

\bibitem{Chiles12}
J.~P. Chil{\`e}s and P.~Delfiner.
\newblock {\em Geostatistics: Modeling Spatial Uncertainty}.
\newblock Wiley, New York, 2nd edition, 2012.

\bibitem{Rasmussen06}
C.~E. Rasmussen and C.~K.~I. Williams.
\newblock {\em Gaussian Processes for Machine Learning}.
\newblock MIT Press, Massachusetts Institute of Technology, 2006.

\bibitem{Christakos92}
G.~Christakos.
\newblock {\em Random Field Models in Earth Sciences}.
\newblock Academic Press, San Diego, 1992.

\bibitem{Cressie11}
N.~Cressie and C.~L. Wikle.
\newblock {\em Statistics for Spatio-temporal Data}.
\newblock John Wiley and Sons, New York, 2011.

\bibitem{Deiaco02}
S.~De~Iaco, D.~E. Myers, and D~Posa.
\newblock Nonseparable space-time covariance models: some parametric families.
\newblock {\em Mathematical Geology}, 34(1):23--42, 2002.

\bibitem{Kolovos04}
A.~Kolovos, G.~Christakos, D.T. Hristopulos, and M.~L. Serre.
\newblock Methods for generating non-separable spatiotemporal covariance models
  with potential environmental applications.
\newblock {\em Advances in Water Resources}, 27(8):815--830, 2004.

\bibitem{dth17b}
E.~A. Varouchakis and D.~T. Hristopulos.
\newblock Comparison of spatiotemporal variogram functions based on a sparse
  dataset of groundwater level variations.
\newblock {\em Spatial Statistics}, 34:100245, 2019.

\bibitem{Gneiting06}
T.~Gneiting, M.~G. Genton, and P.~Guttorp.
\newblock Geostatistical space-time models, stationarity, separability, and
  full symmetry.
\newblock In B.~Finkelst\'{a}dt, L.~Held, and V.~Isham, editors, {\em
  Statistical Methods for Spatio-Temporal Systems}, volume 107, pages 151--175.
  Chapman \& Hall, 2006.

\bibitem{Sun12}
Y.~Sun, B.~Li, and M.~G. Genton.
\newblock Geostatistics for large datasets.
\newblock In E.~Porcu, J.–M. Montero, and M.~Schlather, editors, {\em
  Advances and Challenges in Space-time Modelling of Natural Events}, Lecture
  Notes in Statistics, pages 55--77. Springer Berlin Heidelberg, 2012.

\bibitem{Mussardo10}
G.~Mussardo.
\newblock {\em Statistical Field Theory}.
\newblock Oxford University Press, Oxford, 2010.

\bibitem{Rue05}
H.~Rue and L.~Held.
\newblock {\em Gaussian {M}arkov Random Fields: Theory and Applications}.
\newblock Chapman and Hall/CRC, Boca Raton, FL, 2005.

\bibitem{dth15}
D.~T. Hristopulos.
\newblock Stochastic local interaction {(SLI) model: Bridging} machine learning
  and geostatistics.
\newblock {\em Computers \& Geosciences}, 85(Part B):26--37, 2015.

\bibitem{dth17}
D.~T. Hristopulos and I.~C. Tsantili.
\newblock Space--time covariance functions based on linear response theory and
  the turning bands method.
\newblock {\em Spatial Statistics}, 22, Part 2:321--337, 2017.

\bibitem{koutroulis10}
E.~Koutroulis and D.~Kolokotsa.
\newblock Design optimization of desalination systems power-supplied by {PV and
  W/G} energy sources.
\newblock {\em Desalination}, 258(1-3):171--181, 2010.

\bibitem{Bardossy14}
A.~B{\'a}rdossy and G.~Pegram.
\newblock Infilling missing precipitation records--a comparison of a new
  copula-based method with other techniques.
\newblock {\em Journal of hydrology}, 519, Part A:1162--1170, 2014.

\bibitem{Kardar07}
M.~Kardar.
\newblock {\em Statistical Physics of Fields}.
\newblock Cambridge University Press, 2007.

\bibitem{Ising25}
E.~Ising.
\newblock Beitrag zur theorie des ferromagnetismus.
\newblock {\em Zeitschrift f{\"u}r Physik}, 31(1):253--258, 1925.

\bibitem{Besag74}
J.~Besag.
\newblock Spatial interaction and the statistical analysis of lattice systems.
\newblock {\em Journal of the Royal Statistical Society. Series B
  (Methodological)}, 36(2):192--236, 1974.

\bibitem{Nadaraya64}
E.~A. Nadaraya.
\newblock On estimating regression.
\newblock {\em Theory of Probability and its Applications}, 9(1):141--142,
  1964.

\bibitem{Watson64}
G.~S. Watson.
\newblock Smooth regression analysis.
\newblock {\em Sankhya Ser. A}, 26(1):359--372, 1964.

\bibitem{awr00}
G.~Christakos, D.~T. Hristopulos, and P.~Bogaert.
\newblock On the physical geometry concept at the basis of space/time
  geostatistical hydrology.
\newblock {\em Advances in Water Resources}, 23(8):799--810, 2000.

\bibitem{Christakos17}
G.~Christakos.
\newblock {\em Spatiotemporal Random Fields: Theory and Applications}.
\newblock Elsevier, Amsterdam, Netherlands, 2017.

\bibitem{rf}
M.~Schlather, A.~Malinowski, M.~Oesting, D.~Boecker, K.~Strokorb, S.~Engelke,
  J.~Martini, F.~Ballani, O.~Moreva, J.~Auel, P.~J. Menck, S.~Gross, U.~Ober,
  {C. Berreth}, K.~Burmeister, J.~Manitz, P.~Ribeiro, R.~Singleton, B.~Pfaff,
  and {R Core Team}.
\newblock {\em {RandomFields}: Simulation and Analysis of Random Fields}, 2017.
\newblock R package version 3.1.50.

\bibitem{Schlather15}
M.~Schlather, A.~Malinowski, P.~J. Menck, M.~Oesting, and K.~Strokorb.
\newblock Analysis, simulation and prediction of multivariate random fields
  with package {RandomFields}.
\newblock {\em Journal of Statistical Software}, 63(8):1--25, 2015.

\bibitem{Pebesma16}
E.~Pebesma and G.~Heuvelink.
\newblock Spatio-temporal interpolation using gstat.
\newblock {\em RFID Journal}, 8(1):204--218, 2016.

\bibitem{Wikle19}
C.~K. Wikle, A.~Zammit-Mangion, and N.~Cressie.
\newblock {\em Spatio-temporal Statistics with R}.
\newblock CRC Press, Boca Raton, FL, 2019.

\bibitem{copernicus18}
Copernicus Climate Change~Service C3S.
\newblock {ERA5: Fifth generation of ECMWF} atmospheric reanalyses of the
  global climate, 2018.
\newblock Data retrieved from:
  \url{https://cds.climate.copernicus.eu/cdsapp\#!/home}.

\bibitem{Romero17}
D.~Romero, M.~Ma, and G.~B. Giannakis.
\newblock Kernel-based reconstruction of graph signals.
\newblock {\em IEEE Transactions on Signal Processing}, 65(3):764--778, 2017.

\end{thebibliography}


\end{document}